\documentclass[11pt]{article}
\usepackage[centertags]{amsmath}
\usepackage{amsmath}
\usepackage{amssymb}
\usepackage{amsthm}
\usepackage{color}
\usepackage{latexsym}
\usepackage[T1]{fontenc}
\usepackage[latin1]{inputenc}
\usepackage{fancyhdr}
\usepackage[colorlinks=true,linkcolor=blue,citecolor=blue,urlcolor=blue]{hyperref}

\numberwithin{equation}{section}
\newcommand{\be}{\begin{equation}}
\newcommand{\ee}{\end{equation}}
\newcommand{\dint}{\displaystyle\int}
\newcommand{\1}{1\!\!1}
\newtheorem{definition}{Definition}[section]
\newtheorem{theorem}{Theorem}[section]
\newtheorem{proposition}{Proposition}[section]
\newtheorem{lemma}{Lemma}[section]
\newtheorem{remark}{Remark}[section]

\newcommand{\ind}{\1}

\newcommand{\al}{\alpha}

\def \R{\mathbb{R}}

\def \N{\mathbb{N}}
\def \E{\mathbb{E}}

\def \L{\mathbb{L}}

\def \bf{\textbf}

\begin{document}

\title{One dimensional reflected BSDEs with two barriers under logarithmic growth and applications}
\author{ Brahim EL ASRI \thanks{Universit\'e Ibn Zohr, Equipe. Aide \`a la decision,
		ENSA, B.P.  1136, Agadir, Maroc. e-mail: b.elasri@uiz.ac.ma } \,,\,\, Khalid OUFDIL \thanks{Universit\'e Ibn Zohr, Equipe. Aide \`a la decision,
		ENSA, B.P.  1136, Agadir, Maroc. e-mail: khalid.ofdil@gmail.com}\thanks{K. Oufdil is supported in part by the National Center for Scientific and Technical Research (CNRST), Morocco.}\,\,\,
	\, and \, Nacer OURKIYA \thanks{Universit\'e Ibn Zohr, Equipe. Aide \`a la decision,
		ENSA, B.P.  1136, Agadir, Maroc. e-mail: nacer.ourkiya@edu.uiz.ac.ma}}
\date{}
\maketitle \noindent {\bf{Abstract.}} In this paper we  deal with the problem of the existence and the uniqueness of a solution for one dimensional reflected backward stochastic differential equations with two strictly separated barriers when the generator is allowing a logarithmic growth $(|y||\ln|y||+|z|\sqrt{|\ln|z||})$ in the state variables $y$ and $z$. The terminal value $\xi$  and the obstacle processes $(L_t)_{0\leq t\leq T}$ and $(U_t)_{0\leq t\leq T}$ are $L^p$-integrable for a suitable $p > 2$. The main idea is to use the concept of local solution to construct the global one. As applications, we broaden the class of functions for which mixed zero-sum stochastic differential games admit an optimal strategy and the related double obstacle partial differential equation problem has a unique viscosity solution.\\




\medskip
\section{Introduction}

In this paper we are concerned with the problem of the existence and uniqueness of a solution
for one dimensional reflected backward stochastic differential equations (BSDEs for short) driven by Brownian motion $(B_t)_{t\leq T}$ with two continuous reflecting barriers $L:=(L_t)_{t\leq T}$ and $U:=(U_t)_{t\leq T}$ and whose coefficient and terminal value are $f$ and $\xi$ respectively. Meaning, we want to show the existence of a unique quadruple $(Y,Z,K^+,K^-)$ of $\mathcal{F}_t$-adapted processes such that:

\begin{equation}\label{1}
\left\{\begin{array}{l}
Y_t =\xi + \int_t^T
f(s,Y_s,Z_s) ds +(K^+_T-K^+_t)-(K^-_T-K^-_t)- \int_t^T Z_s dB_s,\qquad
t\in[0,T];\\
\forall t\in[0,T]\qquad L_t\leq Y_t\leq U_t;\\ \int_{0}^{T}\left(Y_s-L_s\right)dK^+_s=\int_{0}^{T}\left(U_s-Y_s\right)dK^-_s=0.
\end{array}\right.
\end{equation}
In the framework of a Brownian filtration, the notion of BSDEs was first introduced by Pardoux and Peng \cite{PP}. Then in \cite{EKPPQ}, El-Karoui et al. introduced BSDEs with a lower obstacle $L:=(L_t)_{t\leq T}$ where the solution $Y$ is assumed to be above $L$, after that, Cvitanic and Karatzas \cite{CK} generalized  these results to BSDEs with two barriers (upper and lower).
Due to their appearance in many finance problems such as the model behind the Black and Scholes formula for the pricing and hedging
of options in mathematical finance, as well as their many applications in other several problems; optimal switching, stochastic games, non-linear PDEs...etc (see \cite{DH,EKPPQ, EPQ, LM} and the references therein).  Many authors have attempted to improve the result of \cite{CK} and establish the existence and the uniqueness of the solution by focusing on weakening the Lipschitz property of the coefficient or the square integrability of the data (see \cite{BHY} for the later).

The main objective of this paper is to show the existence and the uniqueness of the solution for BSDEs with two reflecting barriers with a generator allowing a logarithmic growth in the state variables $y$ and $z$:
$$ \mid f(t,\omega,y,z)\mid \leq |\eta_t|+c_0|y||\ln|y||+c_1|z|\sqrt{\vert\ln(|z|)\vert}~~~~\forall (t,\omega,y,z)\in[0,T]\times\Omega\times\R\times\R^d,$$
with the terminal data $\xi$ and the barriers being merely $p$-integrable (with $ p > 2$). For example, let   $f(y)=-Ky\ln|y|$, and let us consider the following BSDE,
\begin{equation}\label{ly}
	\left\{\begin{array}{l}
		Y_t =\xi + \int_t^T
		f(Y_s) ds +(K^+_T-K^+_t)-(K^-_T-K^-_t)- \int_t^T Z_s dB_s,\qquad
		t\in[0,T];\\
		\forall t\in[0,T]\qquad L_t\leq Y_t\leq U_t;~~ \int_{0}^{T}\left(Y_s-L_s\right)dK^+_s=\int_{0}^{T}\left(U_s-Y_s\right)dK^-_s=0.
	\end{array}\right.
\end{equation}
The generator in \eqref{ly} is not locally monotone nor of sublinear growth in the $y$-variable, moreover, its growth is big power than $y$. The logarithmic
nonlinearity $y \ln|y|$ which appears in \eqref{ly} is interesting in itself and in our knowledge it has not been covered yet, the same thing goes for $f(z)=|z|\sqrt{|\ln|z||}$. As we can see, our assumption covers both cases; $f(y)=-Ky\ln|y|$ and $f(z)=|z|\sqrt{|\ln|z||}$.

Moreover, we also impose an other assumption on $f$ (see, (\textbf{H.4}) below) which is local in $y$, $z$ and also in $\omega$, this enables us to cover certain BSDEs with stochastic monotone generators.

There are mainly two reasons why we study this kind of a problem. The first one is zero-sum games, Dynkin type or of mixed ones, where we broaden the class of data for which those games have a value. It is well known that double barrier reflected BSDEs are connected with mixed zero-sum
games which we describe briefly. Assume that we have a stochastic system whose dynamic $(x_t)_{t\leq T}$ satisfies:
$$x_t=x_0+\int_{0}^{t}\varphi(s,x_s,u_s,v_s)ds+\int_{0}^{t}\sigma(s,x_s)dB_s, ~~t\in[0,T]~~\text{and}~~ x\in\R^d,$$
$\varphi$ is the drift of the system and the stochastic processes $(u_t)_{t\leq T}$ and $(v_t)_{t\leq T}$ are adapted and stand for, respectively, the intervention functions of two agents $A_1$ and $A_2$ on that system; (the system could be for example a stock market and $A_1$ and $A_2$ are two traders). Moreover, the two agents can exit the system whenever they want, meaning, they can stop controlling at stopping times $\tau$ and $\sigma$. However, their actions are not free and their advantages are antagonistic, i.e., there is a payoff $J(u,\tau;v,\sigma)$ between them
 such that
$$J(u,\tau;v,\sigma)=\E^{(u,v)}\left[\int_{0}^{\tau\wedge\sigma}h(s,x,u_s,v_s)ds+L_\sigma\ind_{\{ \sigma\leq \tau<T \}}+U_\tau\ind_{\{ \tau<\sigma\}}+\xi\ind_{\{ \tau\wedge\sigma=T \}}\right],$$
where $h$ is the instantaneous reward of $A_2$, $L$ (resp. $U$) is the reward if $A_2$
decides to stop at $\sigma$ (resp. $\tau$ ) before the terminal time $T$ and $\xi$ is the reward if he decides
to stay until $T$.

The first (resp. second) player chooses a pair $(u, \tau)$ (resp. $(v, \sigma))$ of continuous
control and stopping time, and looks for minimizing (resp. maximizing) this payoff, meaning we aim to find a pair of strategies $(u^*,\tau^*)$ and $(v^*,\sigma^*)$ for $A_1$ and $A_2$ respectively such that $J(u^*,\tau^*;v,\sigma)\leq J(u^*,\tau^*;v^*,\sigma^*)\leq J(u,\tau;v^*,\sigma^*).$ The main idea is we characterize the value function as a solution of a specific reflected BSDE with two barriers.
This problem has already been studied, for example, in \cite{Gam1} when $\sigma^{-1}$,  $\varphi$ and $h$ are  bounded and in \cite{Gam} when $\sigma^{-1}\varphi$ is bounded and $h$ is of linear growth with respect to the $x$-variable. We however, consider the case when $h$ and $\varphi$ are of linear growth with respect to the $x$-variable.

The second reason for considering this problem is to weaken the hypotheses under which the two obstacle parabolic partial differential variational inequality has a unique solution in the viscosity sense.
We consider for example the Markovian of the BSDEs \eqref{ly}, which is defined by the system SDE-BSDE:
\begin{equation}\label{lyy}
\left\{
\begin{array}{l}
X^{t,x}_s=x+\int_{t}^{s}b({u,X^{t,x}_u})du+\int_{t}^{s}\sigma({u,X^{t,x}_u})dB_u,\\
Y^{t,x}_s = g(X^{t,x}_T) -K \int_s^T Y^{t,x}_u\ln|Y^{t,x}_u|
du+\int_{s}^{T}dK^{+,t,x}_u-\int_s^T dK^{-,t,x}_u- \int_s^T Z^{t,x}_u dB_u,\\
\forall s\in[t,T],~~~~h(s,X^{t,x}_s)\leq Y^{t,x}_s\leq h'(s,X^{t,x}_s),\\
\int_{t}^{T} \left( Y^{t,x}_s-h(s,X^{t,x}_s)\right)dK^{+,t,x}_s=\int_{t}^{T} \left(h'(s,X^{t,x}_s)- Y^{t,x}_s\right)dK^{-,t,x}_s=0.
\end{array}\right.
\end{equation}

The system of double obstacle variational inequality associated with \eqref{lyy} is given by
\begin{equation}\label{kk}
	\left\{
	\begin{array}{l}
\min\left[u(t,x)-h(t,x),\max\bigg\{
-\frac{\partial u}{\partial t}(t,x)-\mathcal{L}u(t,x)\right.\\
\left.\qquad\qquad\qquad+Ku(t,x)\ln|u(t,x)|,u(t,x)-h'(t,x)\bigg\}\right]=0,~~ (t,x)\in[0,T)\times \R ^d;\\
u(T,x)=g(x),~~\forall x\in\R ^d,
	\end{array}\right.
\end{equation}
where $$\mathcal{L}=\frac{1}{2}\sum_{i,j=1}^{d}\left(\left(\sigma\sigma^*\right)\left(t,x\right)\right)_{i,j}\frac{\partial^2}{\partial x_i\partial x_j}+\sum_{i=1}^{d}\left(b(t,x)\right)_i\frac{\partial}{\partial x_i}.$$
The logarithmic nonlinearity $u \ln|u|$ is  interesting on its own, since it is neither locally Lipschitz nor uniformly continuous.

This paper is organized as follows. In Section 2, we present the notations and the assumptions used through out the paper. Moreover, we give some preliminary results that would be useful in this paper. In Section 3, we show the existence of a local solution for the two barriers reflected BSDE. Later we show the existence and the uniqueness of the solution for \eqref{1}. In Section 4, we apply the obtained results and we prove that the value function of a mixed zero-sum stochastic differential game problem can be characterized as the solution of a specific BSDE with two barriers.  In Section 5, we show that, provided the problem is formulated within a Markovian framework, the solution of
the reflected BSDE provides a probabilistic representation for the unique viscosity
solution of the related obstacle parabolic partial differential variational inequality.
\section{Notations, Assumptions and Preliminary results}
\subsection{Notations}
Let $(\Omega, \mathcal{F}, P)$ be a fixed probability space on which is defined a standard $d$-dimensional Brownian motion $B=(B_t)_{0\leq t\leq T}$ whose natural filtration is $(\mathcal{F}_t^0:=\sigma \{B_s, s\leq t\})_{0\leq t\leq T}$. Let $\mathcal{F}=(\mathcal{F}_t)_{0\leq t\leq T}$ be
the completed filtration of $(\mathcal{F}_t^0)_{0\leq t\leq T}$ with the $P$-null sets of $\mathcal{F}$. \\

Next for any $p>0$
\begin{itemize}
\item $\mathcal{S}^p$ be the space of $\R $-valued $\mathcal{F}_t$-adapted and continuous processes $\left(Y_t\right)_{t\in[0,T]}$ such that $$||Y||_{\mathcal{S}^p}=\E\left[\sup_{t\leq T}|Y_t|^p\right]^{\frac{1}{p}}<+\infty.$$
\item $\mathcal{M}$ denote the set $\mathcal{P}$-measurable processes $\left(Z_t\right)_{t\in[0,T]}$ with values in $\R^d$  such that:
$$\int_{0}^{T}|Z_s|^2ds<+\infty~~~~\text{P-a.s.};$$
and $\mathcal{M}^p$ is a subset of $\mathcal{M}$ such that:
 $$||Z||_{\mathcal{M}^p}=\E\left[\left(\int_{0}^{T}|Z_s|^2ds\right)^{\frac{p}{2}}\right]^{\frac{1}{p}}<+\infty.$$
\item $\mathcal{A}$ be the set of adapted continuous non decreasing processes $\left(K_t\right)_{t\in[0,T]}$ such that $K_0=0$ and $K_T<+\infty,$ P-a.s. and $\mathcal{A}^p$ is the subset of $\mathcal{A}$ such that $\E\left[K_T^p\right]<+\infty$.
\end{itemize}
\subsection{Assumptions}
Now we are given four data:
\begin{itemize}
\item  $\xi$ is an $\R$-valued and $\mathcal{F}_T$-measurable random variable.
\item  $f:[0,T]\times\Omega\times\R\times\R^d\rightarrow\R$  be a random function which associates $(t,\omega,y,z)$ with $f(t,\omega,y,z)$.
\item  $L:=(L_t)_{0\leq t\leq T}$ and $U:=(U_t)_{0\leq t\leq T}$ are two continuous progressively measurable $\R$-valued processes.
\end{itemize}
On the data $\xi$, $f$, $L$ and $U$ we make the following assumptions:
\begin{enumerate}
\item[(\bf {H.1})] \ \ There exists a positive constant $\lambda$ large enough 
such that $$\E \left[ |\xi|^{e^{\lambda T}+1}\right]<+\infty.$$

\item[(\bf {H.2})]\ \ The two barriers $(L_t)_{0\leq t\leq T}$ and $(U_t)_{0\leq t\leq T}$ satisfy $L_t<U_t,~\forall t\in [0,T]$ and $L_T\leq\xi\leq U_T$. In addition for $p\in]1,2[$ we have $$\E\left[\sup_{0\leq t\leq   T}\left((L^+_t)^{e^{\lambda T}+1}\right)^{\frac{p}{p-1}}\right]<+\infty ~~\text{and}~~ \E\left[\sup_{0\leq t\leq T}\left((U^-_t)^{e^{\lambda T}+1}\right)^{\frac{p}{p-1}}\right]<+\infty,$$
where $L^+=L\vee 0$ and $U^-=(-U)\vee 0.$

\item[(\bf{H.3})]
\begin{enumerate}\item[(i)]  $f$ is continuous in $(y,z)$ for almost all $(t,\omega)$.
    \item[(ii)] There exist three positive constants $c_0$, $\lambda$ (large enough) and $c_1$ and a process $(\eta_t)_{t\leq T} $ such that:
    $$\mid f(t,\omega,y,z)\mid \leq |\eta_t|+c_0|y||\ln|y||+c_1|z|\sqrt{\vert\ln(|z|)\vert}~~~~\forall (t,\omega,y,z)\in[0,T]\times\Omega\times\R\times\R^d, $$
    and $$\E\left[\int_{0}^{T} |\eta_s|^{ e^{\lambda T}+1}ds\right]<+\infty.$$
    \end{enumerate}
\item[(\bf {H.4})]
\par
  \ There exist $ v\in\L^{q'}(\Omega\times [0, T]; \R_
 +)$ (for some $q'>0$), a real valued sequence $(A_N)_{N>1}$ and constants
 $M \in\R_+$, $r>0$ such that:
 \begin{itemize}
 \item[(i)] $\forall N>1$, \quad $1<A_N\leq N^{r}.$
 \item[(ii)] $\lim\limits_{N\rightarrow\infty} A_N =+
 \infty .$
 \item[(iii)] For every $N\in\N,\; \hbox {and every} \ y,\; y',\; z,\; z'
 \;\hbox{such that}\; \mid y\mid,\; \mid y'\mid,\; \mid z\mid, \;\mid
 z'\mid\leq N$, we have:
 \begin{multline*}
  \big(y-y^{\prime}\big)
  \big(f(t,\omega,y,z)-f(t,\omega,y^{\prime},z')\big)
  \1_{\{v_t(\omega)\leq N\}}\\ \leq M\left( \mid
   y-y^{\prime}\mid^{2}\ln A_{N}+ \mid y-y^{\prime}\mid\mid
    z-z^{\prime}\mid\sqrt{\ln A_{N}} +\dfrac{\ln A_{N}}{ A_{N}}\right).
 \end{multline*}
 \end{itemize}
\end{enumerate}

\subsection{Preliminary results}
Now let us define the notion of the local and the global solution of the reflected BSDE associated with the quadruple $(\xi,f,L,U)$ which we consider throughout this paper. We start with the global solution.
\begin{definition}\label{def}
We say that $\{(Y_t,Z_t,K^+_t,K^-_t);0\leq t\leq T\}$ is a solution of the reflected BSDE associated with two continuous barriers $L$ and $U$, a terminal condition
$\xi$ and  a generator $f$ if the followings hold:
\begin{equation}\label{equ}
\left\{\begin{array}{l}
Y\in \mathcal{S}^{e^{\lambda T}+1},~~Z\in\mathcal{M},~~K^\pm\in\mathcal{A};\\
Y_t =\xi + \int_t^T
f(s,Y_s,Z_s) ds +(K^+_T-K^+_t)-(K^-_T-K^-_t)- \int_t^T Z_s dB_s,\qquad
t\in[0,T];\\
\forall t\in[0,T]\qquad L_t\leq Y_t\leq U_t;\\ \int_{0}^{T}\left(Y_s-L_s\right)dK^+_s=\int_{0}^{T}\left(U_s-Y_s\right)dK^-_s=0.
\end{array}\right.
\end{equation}
\end{definition}
Since in many applications, especially in stochastic games or mathematical finance, we don't need strong
integrability conditions on $Z$ and $K^\pm$, that is why we don't require them as one can notice from the definition \ref{def}.\\

Now we define the local solution. In the following $p\in]1,2[$
\begin{definition}\label{diff}
Let $\tau$ and $\gamma$ be two stopping times such that $\tau\leq\gamma$ $P$-a.s. We say that $(Y_t,Z_t,K^+_t,K^-_t)_{0\leq t\leq T}$ is a local solution on $[\tau,\gamma]$ of the reflected BSDE associated with two continuous barriers $L$ and $U$, a terminal condition $\xi$ and a generator $f$ if the followings hold:
\begin{equation}\label{lequ}
\left\{\begin{array}{l}
Y\in \mathcal{S}^{e^{\lambda T}+1},~~Z\in\mathcal{M}^2,~~K^\pm\in\mathcal{A}^p ;\\
Y_t =Y_{\gamma} + \int_t^\gamma
f(s,Y_s,Z_s) ds +(K^+_\gamma-K^+_t)-(K^-_\gamma-K^-_t)- \int_t^\gamma Z_s dB_s,\qquad
\forall t\in[\tau,\gamma];\\
Y_T=\xi;\\
L_t\leq Y_t\leq U_t,\qquad \forall t\in[\tau,\gamma]~~~and~~~\int_{\tau}^{\gamma}\left(Y_s-L_s\right)dK^+_s=\int_{\tau}^{\gamma}\left(U_s-Y_s\right)dK^-_s=0.
\end{array}\right.
\end{equation}
\end{definition}
We first begin with an estimation of $f$ which can easily be proved.
\begin{lemma}\label{estf} If \bf{(H.3)}
 holds, then for any $\alpha\in]1,2[$
\begin{align*}
&  \E \left[ \dint_0^T   | f(s,Y_s,Z_s)|^{\frac{2}{\alpha}} ds \right] \; \leq \;
K\E\left[\dint_0^T
 |{\eta}_s|^2 ds+\sup_{s\leq T}|Y_s|^{\frac{4}{{\alpha}}}
+\dint_0^T\vert Z_s\vert^2 ds\right],
\end{align*}
where $K$ is a positive constant that depends on $c_0$ and $T$.
\end{lemma}
We now introduce the comparison result established in [\cite{EO}, Theorem 4.1] which also holds in our setting.
\begin{proposition}\label{comp}
	Let $(\xi,f,L)$ and $(\xi',f',L')$ be two sets of data that satisfies all the assumptions;
	(\textbf{H.1}), (\textbf{H.2}), (\textbf{H.3}), and (\textbf{H.4}). And suppose in addition the followings:
	\begin{itemize}
		\item[(i)] $\xi\leq \xi'$ $P$-a.s.
		\item[(ii)] $f(t,y,z)\leq f'(t,y,z)$
		$dP\times dt$ a.e.,
		$\forall (t,y,z)\in[0,T]\times \R\times \R^d. $
		\item[(iii)] $L_t\leq L'_t;\qquad \forall t\in[0,T]$ $P$-a.s.
	\end{itemize}
	Let $(Y,Z,K^+)$  be the  solution of the reflected BSDE with one lower barrier associated with $(\xi,f,L)$ i.e.
	\begin{equation*}
		\left\{\begin{array}{l}
			Y_t =\xi + \int_t^{T}f(s,Y_s,Z_s) ds +K^+_T-K^+_t- \int_t^{T} Z_s dB_s,\qquad
			t\in[0,T];\\
			\forall t\in[0,T]\qquad L_t\leq Y_t;\\ \int_{0}^{T}\left(Y_s-L_s\right)dK^+_s=0,
		\end{array}\right.
	\end{equation*}
	and $(Y',Z',K'^+)$ the solution of the reflected BSDE with one lower barrier associated with $(\xi',f',L')$. Then,
	$$Y_t\leq Y_t',~~~~~~~0\leq t\leq T~~P\text{-a.s.}$$
\end{proposition}
\newpage
\begin{remark}\label{comprem}${}$
\begin{itemize}
\item The comparison result also holds for reflected BSDEs with the one upper barrier $U$, that is, if $(\xi,f,U)$ and $(\xi',f',U')$ are two sets of data satisfying (\textbf{H.1})-(\textbf{H.4}), if moreover $\xi\leq \xi'$, $f(t,y,z)\leq f'(t,y,z)$ and $U\leq U'$, then $P$-a.s.,  $Y_t\leq Y_t',\forall~0\leq t\leq T$, where $(Y,Z,K^-)$ is the solution of the one upper barrier reflected BSDEs associated with $(\xi,f,U)$ i.e.
	\begin{equation*}
	\left\{\begin{array}{l}
		Y_t =\xi + \int_t^{T}f(s,Y_s,Z_s) ds -K^-_T+K^-_t- \int_t^{T} Z_s dB_s,\qquad
		t\in[0,T];\\
		\forall t\in[0,T]\qquad U_t\geq Y_t;\\ \int_{0}^{T}\left(U_s-Y_s\right)dK^-_s=0,
	\end{array}\right.
\end{equation*}

 and $(Y',Z',K'^-)$ is the solution of the one upper barrier reflected BSDEs associated with $(\xi',f',U')$.
\item  If $L=-\infty$, then $K^+=0$ and the comparison theorem holds in the case with no barrier.
\end{itemize}	
\end{remark}
\section{Existence and uniqueness of the solution}
In this section we are going to show the existence and the uniqueness of the solution for \eqref{equ}, but first we show that it has a local solution in the sens of Definition \ref{diff},  before we show that this solution is in fact a global one when the barriers are completely. The main difficulty, in this section, is to show that the solution of the one barrier reflected BSDE studied in \cite{EO} can be obtained using penalization method and comparison theorem, since the authors of \cite{EO} used localization technique to get that result. Actually, we have the following theorem.
\begin{theorem}\label{penloga}
$\forall n\geq 0$, let $(y^n_t,z^n_t)_{t\leq T}$ be the unique solution of the BSDE
\begin{equation}\label{conse}
y^n_t =\xi + \int_t^{T}
\left(f(s,y^n_s,z^n_s)+n(L_s-y^n_s)^+\right) ds - \int_t^{T} z^n_s dB_s,\qquad
t\in[0,T],
\end{equation}
which exists due to [\cite{BKK}, Theorem 2.1]. Then, the processes $(y^n_s, z^n_s, \int_{0}^{s}n(L_r-y^n_r)^+dr)_{s\leq T}$ converges to $(y_s,z_s,k_s)_{s\leq T}$ solution of
\begin{equation}\label{solref}
\left\{\begin{array}{l}
\E\left[\sup\limits_{0\leq s\leq T }|y_s|^{e^{\lambda T}+1}+\int_{0}^{T}|z_s|^2 ds+k^p_T\right]<+\infty,\\
y_t =\xi + \int_t^{T}f(s,y_s,z_s) ds +k_T-k_t- \int_t^{T} z_s dB_s,\qquad
t\in[0,T].\\
\forall t\in[0,T]\qquad L_t\leq y_t;\\ \int_{0}^{T}\left(y_s-L_s\right)dk_s=0.
\end{array}\right.
\end{equation}
\end{theorem}
\noindent
\textbf{Proof.} The first part of Equation \eqref{solref} is a direct consequence of [\cite{EO}, Proposition 3.1]; that is,
there
exists a positive constant $C(\lambda,p,c_0,c_1,T)$ such that: $\forall p\in]1,2[$,
\begin{equation}
	\begin{array}{ll}\label{estKYZ}
		\E \left[\sup\limits_{t \in [0,T]} |y_t|^{e^{\lambda t}+1}+ \int_0^T |z_s|^2 ds + k_T^{p}\right] \\ \qquad\quad\leq C(\lambda,p,c_0,c_1,,T)\E \left[1+|\xi|^{e^{\lambda T}+1} +\int_0^T\mid \eta_s \mid^{e^{\lambda s}+1}ds+ \sup\limits_{0\leq t\leq T}\left((L^+_t)^{e^{\lambda t}}\right)^{\frac{p}{p-1}}\right].
	\end{array}
\end{equation}
Next, we define
$$k^n_t=n\int_{0}^{t}(L_s-y^n_s)^+ds, \qquad t\in[0,T].$$ Hence, from \eqref{estKYZ}, we have that for $p\in]1,2[$
\begin{equation}\label{fat}
\E\left[\sup\limits_{0\leq s\leq T}|y^n_s|^{e^{\lambda T}+1}+\int_{0}^{T}|z^n_s|^2 ds+(k^n_T)^p\right]<+\infty,~~~\forall n\geq 0.	
\end{equation}
Note that if we define $f_n(t,y,z)=f(t,x,y)+n(L_t-y)^+$, then $f_n(t,y,z)\leq f_{n+1}(t,y,z)$. Using the comparison Theorem in \cite{EO}, it follows that $y^n_t\leq y^{n+1}_t,$ $0\leq t\leq T$, a.s.
Therefore, by dominated convergence we have
\begin{eqnarray}\label{domconv}
\E\left[\int_{0}^{T}(y_t-y^n_t)^{e^{\lambda T}+1}dt\right]\rightarrow 0 \text{ as } n\rightarrow+\infty.
\end{eqnarray}
The rest of the proof will be divided into two steps.\\

\noindent
\textbf{Step 1.}
We will show that for $p\in]\frac{e^{{\lambda T}+1}}{e^{{\lambda T}+1}-1},2[$
\begin{equation}\label{leem}
\E\left[ \sup\limits_{0\leq s\leq T}\left( \big(
L_s-y_{s}^{n}\big)^+\right)^{\frac{p}{p-1}}\right]^{\frac{p-1}{p}}\rightarrow 0 \text{ as } n\rightarrow +\infty.
\end{equation}
For any $n\geq 0$ and $t\leq T$, we have
\begin{equation}\label{forw}
y^n_t= \xi + \int_{t}^{T}f(s,y^n_s,z^n_s)ds+k^n_t-\int_{t}^{T}z^n_sdB_s.	
\end{equation}
Putting  $g^n_s=f(s,y^n_s,z^n_s)$ and writing \eqref{forw} forwardly we get
$$ k^n_t = y^n_0-y^n_t-\int_{0}^{t}g^n_sds+\int_{0}^{t}z^n_sdB_s.$$
Since from Lemma \ref{estf} and \eqref{fat} for any $n\geq 0$
$$\E\left[\sup\limits_{0\leq t\leq T}|y^n_t|^{e^{\lambda T}+1}+\int_{0}^{T}|g^n_s|^{\frac{2}{\al}}ds+\int_{0}^{T}|z^n_s|^2ds\right]<+\infty,$$
then there exist subsequences and processes $(g_t)_{0\leq t\leq T}$ and $(z_t)_{0\leq t\leq T}$ which are the weak limit of $(g^n_t)_{0\leq t\leq T}$ and $(z^n_t)_{0\leq t\leq T}$ respectively. Henceforth, for any stopping time $\bar{\tau} \leq T$ the following weak convergence holds
$$\int_{0}^{\bar{\tau}}z^n_sdB_s\rightarrow \int_{0}^{\bar{\tau}}z_sdB_s~~\text{and}~~ \int_{0}^{\bar{\tau}}g^n_sds\rightarrow \int_{0}^{\bar{\tau}}g_sds.$$
It follows that
$$ k^n_{\bar{\tau}}\rightarrow k_{\bar{\tau}}= y_0-y_{\bar{\tau}}-\int_{0}^{{\bar{\tau}}}g_sds+\int_{0}^{{\bar{\tau}}}z_sdB_s.$$
Now for any stopping times $\bar{\sigma}\leq \bar{\tau}\leq T$ we have $k^n_{\bar{\sigma}}\leq k^n_{\bar{\tau}}$, therefore,  it holds true that $k_{\bar{\sigma}}\leq k_{\bar{\tau}}$. Hence,  $(k_t)_{0\leq t\leq T}$ is an increasing process. Additionally, $\E\left[(k_T)^p\right]\leq \underset{n\rightarrow+\infty}{\lim\inf}~\E\left[(k^n_T)^p\right]<+\infty$. Henceforth, thanks to the monotonic limit of Peng [\cite{PE}, Lemma 2.2], the processes $(y_t)_{0\leq t\leq T}$ and $(k_t) _{0\leq t\leq T}$ are RCLL.

Next, due to the fact that $\E\left[(k^n_T)^p\right]<+\infty$ for any $n\geq 0$, we deduce, in taking the limit $n\rightarrow+\infty$, that
$$\E\left[\int_{0}^{T}(L_s-y_s)^+ds\right]=0.$$
Therefore, $P$-a.s. $y_t\geq L_t$ for any $t<T$. But $\xi\geq L_T$, then $y\geq L$. Hence, $(L_t-y^n_t)^+\downarrow 0$, $0\leq t\leq T$, a.s. and from Dini's theorem the convergence is uniform in $t$. Since $(L_t-y^n_t)^+\leq |L_t|+|y^0_t|$, the result follows.\\

\noindent
\textbf{Step 2.} We will show that $(y^n,z^n,k^n)$ converges to $(y,z,k)$ solution of \eqref{solref}.\\

\noindent
Let $0\leq T'\leq T$ and put $ \Delta_{t}:=\left|
y_{t}^{n}-y_{t}^{m}\right| ^{2}+ (A_N)^{-1} $ and $\Phi(s)=|y_{s}^{n}| + |y_{s}^{m}|+
|z_{s}^{n}| + |z_{s}^{m}|+v_s$. Then, for $C>0$ and  $1<\beta < \min
\{(3-\alpha),2\}$

\begin{eqnarray*}
&& e^{Ct}\Delta_{t}^{\beta\over 2}
+C\int_{t}^{T'}e^{Cs}\Delta_{s}^{\frac{\beta}{2}}ds =
e^{CT'}\Delta_{T'}^{\beta\over 2}\\
&&\qquad\qquad\qquad	+\beta\int_{t}^{T'}e^{Cs}\Delta_{s}^{\frac{\beta}{2}-1} \big(
y_{s}^{n}-y_{s}^{m}\big)
\big(f(s,y_{s}^{n},z_{s}^{n})-
f(s,y_{s}^{m},z_{s}^{m})\big)\ind_{\Phi(s)>N}
ds\\
&&\qquad\qquad\qquad+\beta\int_{t}^{T'}e^{Cs}\Delta_{s}^{\frac{\beta}{2}-1} \big(
y_{s}^{n}-y_{s}^{m}\big)
\big(f(s,y_{s}^{n},z_{s}^{n})-
f(s,y_{s}^{m},z_{s}^{m})\big)\ind_{\Phi(s)\leq N}
ds\\
&&\qquad\qquad\qquad -\frac{\beta}{2}\int_{t}^{T'}e^{Cs}\Delta_{s}^{\frac{\beta}{2}-1}\left|
z_{s}^{n}-z_{s}^{m}\right|^{2}ds
-\beta\int_{t}^{T'}e^{Cs}\Delta_{s}^{\frac{\beta}{2}-1}\left(
y_{s}^{n}-y_{s}^{m}\right) \left(
z_{s}^{n}-z_{s}^{m}\right)dB_{s}\\
&& \qquad\qquad\qquad-\beta(\frac{\beta-2}{2})\int_{t}^{T'}e^{Cs}
\Delta_{s}^{\frac{\beta}{2}-2}
\left((y_{s}^{n}-y_{s}^{m})(z_{s}^{n}-
z_{s}^{m})\right)^{2}ds\\
&&\qquad\qquad\qquad+\beta\int_{t}^{T'}e^{Cs}\Delta_{s}^{\frac{\beta}{2}-1} \big(
y_{s}^{n}-y_{s}^{m}\big)\big(dk^{n}_s-dk^{m}_s \big).
\end{eqnarray*}
First let us deal with
$$\beta\int_{t}^{T'}e^{Cs}\Delta_{s}^{\frac{\beta}{2}-1} \big(
y_{s}^{n}-y_{s}^{m}\big)\big(dk^{n}_s-dk^{m}_s \big).$$
Actually,
\begin{eqnarray*}
&& \beta\int_{t}^{T'}e^{Cs}\Delta_{s}^{\frac{\beta}{2}-1} \big(
y_{s}^{n}-y_{s}^{m}\big)\big(dk^{n}_s-dk^{m}_s \big)\\
&&\qquad\qquad\qquad\qquad = \beta\int_{t}^{T'}e^{Cs}\left(\left|
y_{t}^{n}-y_{t}^{m}\right| ^{2}+ (A_N)^{-1}\right)^{\frac{\beta}{2}-1} \big(
y_{s}^{n}-y_{s}^{m}\big)dk^{n}_s\\
&& \qquad\qquad\qquad\qquad+ \beta\int_{t}^{T'}e^{Cs}\left(\left|
y_{t}^{m}-y_{t}^{n}\right| ^{2}+ (A_N)^{-1}\right)^{\frac{\beta}{2}-1} \big(
y_{s}^{m}-y_{s}^{n}\big)dk^{m}_s
\end{eqnarray*}
Since $dk^n_s=\ind_{\{y^n_s\leq L_s\}}dk^n_s$ and $dk^m_s=\ind_{\{y^m_s\leq L_s\}}dk^m_s$ and the function $x\mapsto\beta e^{Cs}(|x-y|^2+(A_N)^{-1})^{\frac{\beta}{2}-1} \big(x-y\big)$ is non-decreasing, it follows that
 \begin{eqnarray*}
 	&& \beta\int_{t}^{T'}e^{Cs}\Delta_{s}^{\frac{\beta}{2}-1} \big(
 	y_{s}^{n}-y_{s}^{m}\big)\big(dk^{n}_s-dk^{m}_s \big)\\
 	&&\qquad\qquad\qquad\qquad = \beta\int_{t}^{T'}e^{Cs}\left(\left|
 	L_s-y_{s}^{m}\right| ^{2}+ (A_N)^{-1}\right)^{\frac{\beta}{2}-1} \big(
 	L_s-y_{s}^{m}\big)dk^{n}_s\\
 	&& \qquad\qquad\qquad\qquad+ \beta\int_{t}^{T'}e^{Cs}\left(\left|
 	L_s-y_{s}^{n}\right| ^{2}+ (A_N)^{-1}\right)^{\frac{\beta}{2}-1} \big(
 	L_s-y_{s}^{n}\big)dk^{m}_s\\
 	&&\qquad\qquad\qquad\qquad \leq 2\beta e^{CT'} \sup\limits_{0\leq s\leq T}\left(\left(\left|
 	L_s-y_{s}^{m}\right| ^{2}+ (A_N)^{-1}\right)^{\frac{\beta}{2}-1} \big(
 	L_s-y_{s}^{m}\big)^+\right) k^n_T\\
 	&&\qquad\qquad\qquad\qquad + 2\beta e^{CT'} \sup\limits_{0\leq s\leq T}\left(\left(\left|
 	L_s-y_{s}^{n}\right| ^{2}+ (A_N)^{-1}\right)^{\frac{\beta}{2}-1} \big(
 	L_s-y_{s}^{n}\big)^+\right) k^m_T.
 \end{eqnarray*}
Since $\frac{\beta}{2}-1 <0$ and since $\forall t\in [0,T]$
$$ (A_N)^{-1} \leq \left|
L_{t}-y_{t}^{m}\right| ^{2}+ (A_N)^{-1} \text{ and } (A_N)^{-1} \leq \left|
L_{t}-y_{t}^{n}\right| ^{2}+ (A_N)^{-1}$$
then $$\left(\left|
L_{t}-y_{t}^{n}\right| ^{2}+ (A_N)^{-1}\right)^{\frac{\beta}{2}-1}\leq (A_N)^{\frac{2-\beta}{2}} \text{ and } \left(\left|
L_{t}-y_{t}^{m}\right| ^{2}+ (A_N)^{-1}\right)^{\frac{\beta}{2}-1}\leq (A_N)^{\frac{2-\beta}{2}}.$$
It follows that
\begin{eqnarray*}
	&&  \beta\int_{t}^{T'}e^{Cs}\Delta_{s}^{\frac{\beta}{2}-1} \big(
	y_{s}^{n}-y_{s}^{m}\big)\big(dk^{n}_s-dk^{m}_s \big)\\
	&&\qquad\qquad\qquad\leq 2(A_N)^{\frac{2-\beta}{2}}\beta e^{CT'} \sup\limits_{0\leq s\leq T} \big(
	L_s-y_{s}^{m}\big)^+ k^n_T\\
	&&\qquad\qquad\qquad + 2(A_N)^{\frac{2-\beta}{2}}\beta e^{CT'} \sup\limits_{0\leq s\leq T} \big(
	L_s-y_{s}^{n}\big)^+ k^m_T.
\end{eqnarray*}
Next we put
\begin{eqnarray*}
&& J_1 = \beta\int_{t}^{T'}e^{Cs}\Delta_{s}^{\frac{\beta}{2}-1} \big(
y_{s}^{n}-y_{s}^{m}\big)
\big(f(s,y_{s}^{n},z_{s}^{n})-
f(s,y_{s}^{m},z_{s}^{m})\big)\ind_{\Phi(s)>N}
ds\\
&& J_2 = \beta\int_{t}^{T'}e^{Cs}\Delta_{s}^{\frac{\beta}{2}-1} \big(
y_{s}^{n}-y_{s}^{m}\big)
\big(f(s,y_{s}^{n},z_{s}^{n})-
f(s,y_{s}^{m},z_{s}^{m})\big)\ind_{\Phi(s)\leq N}
ds.
\end{eqnarray*}
Let $\kappa
=3-\alpha-\beta$. Since $\frac
{(\beta-1)}{2}+\frac{\kappa}{2}+\frac{\alpha}{2}=1$, we
use H\" older's inequality to obtain
\begin{eqnarray*}
&& J_1\leq \beta e^{CT'} \dfrac{1}{N^\kappa}
\left[\int_{t}^{T'} \Delta_{s} ds \right]^{\frac{\beta-1}{2}}
\times\left[\int_{t}^{T'}{\Phi(s)}^2 ds \right]^{\frac{\kappa}{2}} \times\left[\int_{t}^{T'}|f(s,y_{s}^{n},z_{s}^{n})-f(s,y_{s}^{m},z_{s}^{m})|^{2\over \alpha}
ds \right]^{\frac{\alpha}{2}}.
\end{eqnarray*}
For $J_2$ we use assumption \textbf{(H.4)} and we obtain
\begin{eqnarray*}
&&J_{2} \leq \beta M
	\int_{t}^{T'}e^{Cs}\Delta_{s}^{\frac{\beta}{2}-1}
	 \bigg[|y^n_s-y^m_s|^2\ln
	A_{N}+ \frac{\ln A_N}{A_N}+
	|y^n_s-y^m_s||z^n_s-z^m_s|\sqrt{\ln
		A_{N}} \bigg]\1_{\{\Phi(s)\leq N\}}ds
	\\ &&\qquad \leq
	\beta M\int_{t}^{T'}e^{Cs}\Delta_{s}^{\frac{\beta}{2}-1}
	 \bigg[\Delta_{s}\ln A_{N}+|y^n_s-y^m_s||z^n_s-z^m_s|\sqrt{\ln
		A_{N}}  \bigg]\1_{\{\Phi(s)\leq N\}}ds.
	\end{eqnarray*}
We apply Lemma 4.6 in \cite{BE}  and we choose $C=C_N=\dfrac{2M^2\beta}{\beta -1} \ln A_{N}$ to get
\begin{eqnarray*}
	&& e^{C_Nt}\Delta_{t}^{\beta\over 2} + \dfrac{\beta(\beta-1)}{4}
	\int_{t}^{T'}e^{C_Ns}\Delta_{s}^{\frac{\beta}{2}-1}\left|
	z_{s}^{n}-z_{s}^{m}\right|^{2}ds\\
	&&\qquad\qquad\qquad\leq e^{C_NT'}\Delta_{T'}^{\beta\over 2} -\beta\int_{t}^{T'}
	e^{C_Ns}\Delta_{s}^{\frac{\beta}{2}-1}\left(
	y_{s}^{n}-y_{s}^{m}\right) \left(
	z_{s}^{n}-z_{s}^{m}\right)dB_{s}\\
	&&\qquad\qquad\qquad+\beta e^{C_NT'} \dfrac{1}{N^\kappa}
	\left[\int_{t}^{T'} \Delta_{s} ds \right]^{\frac{\beta-1}{2}}
	\times\left[\int_{t}^{T'}{\Phi(s)}^2 ds \right]^{\frac{\kappa}{2}}\\
	&&\qquad\qquad\qquad \times\left[\int_{t}^{T'}|f(s,y_{s}^{n},z_{s}^{n})-f(s,y_{s}^{m},z_{s}^{m})|^{2\over \alpha}
	ds \right]^{\frac{\alpha}{2}}\\
&&\qquad\qquad\qquad+ 2(A_N)^{\frac{2-\beta}{2}}\beta e^{CT'} \sup\limits_{0\leq s\leq T} \big(
L_s-y_{s}^{m}\big)^+ k^n_T\\
&&\qquad\qquad\qquad + 2(A_N)^{\frac{2-\beta}{2}}\beta e^{CT'} \sup\limits_{0\leq s\leq T} \big(
L_s-y_{s}^{n}\big)^+ k^m_T.
\end{eqnarray*}
Therefore, from Burkholder's inequality, H\"older's inequality, \eqref{fat} and Lemma \ref{estf} there exists a universal constant $\ell$ such that for $p\in]\frac{e^{{\lambda T}+1}}{e^{{\lambda T}+1}-1},2[$

\begin{eqnarray*}
&& \E \left[\sup_{(T'-\delta')^{+}\leq t \leq T'}\vert
y_{t}^{n}-y_{t}^{m}\vert^\beta\right] + \E\left[
\int_{(T'-\delta')^{+}}^{T'}\dfrac{\left|
	z_{s}^{n}-z_{s}^{m}\right|^{2}}{\left(\vert
	y_{s}^{n}-y_{s}^{m}\vert^{2}+ \nu_{R}
	\right)^{\frac{2-\beta}{2}}}ds\right]\\
&&\qquad\qquad\qquad\leq \ell\left( e^{C_N\delta'}\E\left[ \vert
y_{T'}^{n}-y_{T'}^{m}\vert^\beta\right] + \dfrac{A_N^{\frac{2M^2\delta'\beta}{\beta-1}}}{(A_N)^{\frac{\beta}{2}}}+ \dfrac{A_N^{\frac{2M^2\delta'\beta}{\beta-1}}}
{(A_N)^{\frac{\kappa}{r}}}  \right)\\
&&\qquad\qquad\qquad+ \ell (A_N)^{\frac{2-\beta}{2}}e^{C_N\delta'}\E\left[ \sup\limits_{0\leq s\leq T}\left( \big(
L_s-y_{s}^{n}\big)^+\right)^{\frac{p}{p-1}}\right]^{\frac{p-1}{p}}\\
&&\qquad\qquad\qquad+ \ell (A_N)^{\frac{2-\beta}{2}}e^{C_N\delta'}\E\left[ \sup\limits_{0\leq s\leq T}\left( \big(
L_s-y_{s}^{m}\big)^+\right)^{\frac{p}{p-1}}\right]^{\frac{p-1}{p}}
\end{eqnarray*}
where $\nu_{R} = \sup \left\{(A_N)^{-1}, N\geq R
\right\}$.  Hence for $ \delta' < (\beta
 -1)\min\left(\frac{1}{4M^2},\frac{\kappa}{2rM^2\beta}\right)$
 we derive
 \begin{equation*}
 	\lim_{N\rightarrow+\infty}\dfrac{A_N^{\frac{2M^2\delta'\beta}{\beta-1}}}{(A_N)^{\frac{\beta}{2}}}= 0\qquad\qquad \text{and} \qquad\qquad \lim_{N\rightarrow+\infty}\displaystyle\dfrac{A_N^{\frac{2M^2\delta'\beta}{\beta-1}}}
 	{(A_N)^{\frac{\kappa}{r}}}
 	\displaystyle=0.
 \end{equation*}
It follows then from \eqref{leem} that,
\begin{equation}\label{key}
\limsup_{n,m\rightarrow +\infty} \E \left[\sup_{(T'-\delta')^{+}\leq t\leq T'}\vert
y_{t}^{n}-y_{t}^{m}\vert^\beta\right] \leq \varepsilon +\ell
e^{C_N\delta'} \limsup_{n,m\rightarrow +\infty} \E \left[\vert
y_{T'}^{n}-y_{T'}^{m}\vert^\beta\right].	
\end{equation}
Taking successively $T'=T$,
$T'=(T-\delta')^+$, $T'=(T-2\delta')^{+}...$ in \eqref{key} we get
\begin{equation}\label{beta}
\lim_{n,m\rightarrow +\infty} \E \left[\sup_{0\leq t \leq T}\vert
y_{t}^{n}-y_{t}^{m}\vert^\beta\right] =0.
\end{equation}
Next, we prove that
\begin{equation}\label{exzp}
	\lim\limits_{n,m\rightarrow+\infty}\E\left[\int_{0}^{T}\left| z_{s}^{n}-z_{s}^{m}\right|
	^{2}ds\right]=0.
\end{equation}
It follows from It\^{o}'s formula that:
\begin{align}\label{Exz}
	&\left| y_{0}^{n}-y_{0}^{m}\right|
	^{2}+\int_{0}^{T}\left| z_{s}^{n}-z_{s}^{m}\right|
	^{2}ds
	\\\nonumber&\qquad = 2\int_{0}^{T}\big( y_{s}^{n}-y_{s}^{m}\big)
	\big(f(s,y_{s}^{n},z_{s}^{n})-
	f(s,y_{s}^{m},z_{s}^{m})\big)
	ds
	\\\nonumber&\qquad+2\int_{0}^{T}\big( y_{s}^{n}-y_{s}^{m}\big)\big(dk^{n}_s-dk^{m}_s \big)-2\int_{0}^{T}\left(y_{s}^{n}-y_{s}^{m}\right) \left(
	z_{s}^{n}-z_{s}^{m}\right)dB_{s}.
\end{align}
First we argue that the third term of the right side in \eqref{Exz} is a martingale. We can deduce from Burkholder-Davis-Gundy's inequality and Lemma \ref{fat} that there exists a constant $c>0$ such that:
\begin{eqnarray}\label{mar}
	&& \E\left[\sup_{0\leq t\leq T}\left|\int_{0}^{t}\left(y_{s}^{n}-y_{s}^{m}\right)\left(z_{s}^{n}-z_{s}^{m}\right)dB_{s}\right|\right]\\ \nonumber
	&&\qquad\qquad\qquad\qquad\qquad\leq c\E \left[\sup_{0\leq s\leq T}|y_{s}^{n}-y_{s}^{m}|^2\right]+c\E\left[\int_{0}^{T}|z_{s}^{n}-z_{s}^{m}|^2ds\right]\\ \nonumber
	&& \qquad\qquad\qquad\qquad\qquad<+\infty.
\end{eqnarray}
Now we deal with the term
$\int_{0}^{T}\big( y_{s}^{n}-y_{s}^{m}\big)\big(dk^{n}_s-dk^{m}_s)$. Actually, since $dk^n_s=\ind_{\{y^n_s\leq L_s\}}dk^n_s$ and $dk^m_s=\ind_{\{y^m_s\leq L_s\}}dk^m_s$ we obtain
\begin{eqnarray}\label{k}
	&& \int_{0}^{T}\big( y_{s}^{n}-y_{s}^{m}\big)\big(dk^{n}_s-dk^{m}_s)\\\nonumber
	&& \qquad\qquad\qquad=\int_{0}^{T}\big( y_{s}^{n}-y_{s}^{m}\big)dk^{n}_s+\int_{0}^{T}\big( y_{s}^{m}-y_{s}^{n}\big)dk^{m}_s\\\nonumber
	&& \qquad\qquad\qquad\leq\int_{0}^{T}\big( L_{s}-y_{s}^{m}\big)dk^{n}_s+\int_{0}^{T}\big( L_{s}-y_{s}^{n}\big)dk^{m}_s\\\nonumber
	&& \qquad\qquad\qquad\leq\int_{0}^{T}\big( L_{s}-y_{s}^{m}\big)^+dk^{n}_s+\int_{0}^{T}\big( L_{s}-y_{s}^{n}\big)^+dk^{m}_s.
\end{eqnarray}
Combining \eqref{Exz}, \eqref{mar} and \eqref{k} to obtain that there exists a constant $c$ such that:
\begin{eqnarray}\label{z0}
	&& \E \left[\int_{0}^{T}\left| z_{s}^{n}-z_{s}^{m}\right|
	^{2}ds\right]\leq c\E\left[\int_{0}^{T}\big| y_{s}^{n}-y_{s}^{m}\big|\big|f(s,y_{s}^{n},z_{s}^{n})-
	f(s,y_{s}^{m},z_{s}^{m})\big|ds \right]\\ \nonumber
	&& \qquad\qquad\qquad\qquad\qquad+c\E\left[\int_{0}^{T}\big( L_{s}-y_{s}^{m}\big)^+dk^{n}_s+\int_{0}^{T}\big( L_{s}-y_{s}^{n}\big)^+dk^{m}_s\right].
\end{eqnarray}
Next by  H\" older's inequality we have
\begin{eqnarray}
	&& \E\left[\int_{0}^{T}\big| y_{s}^{n}-y_{s}^{m}\big|\big|f(s,y_{s}^{n},z_{s}^{n})-
	f(s,y_{s}^{m},z_{s}^{m})\big|ds \right]\\ \nonumber
	 &&\quad\quad\qquad\leq\E\left[\left(\int_{0}^{T}|y_s^n-y_s^m|^{\frac{2}{2-\alpha}}ds\right)^{\frac{2-\alpha}{2}}\left(\int_{0}^{T}\big|f(s,y_{s}^{n},z_{s}^{n})-
	f(s,y_{s}^{m},z_{s}^{m})\big|^{\frac{2}{\alpha}}ds\right)^{\frac{\alpha}{2}}\right]\\ \nonumber
	&&\quad\quad\qquad\leq\E\left[\int_{0}^{T}|y_s^n-y_s^m|^{\frac{2}{2-\alpha}}ds\right]^{\frac{2-\alpha}{2}}\E\left[\int_{0}^{T}\big|f(s,y_{s}^{n},z_{s}^{n})-
f(s,y_{s}^{m},z_{s}^{m})\big|^{\frac{2}{\alpha}}ds\right]^{\frac{\alpha}{2}}.
\end{eqnarray}
We plug the last inequality in \eqref{z0} and we get
\begin{eqnarray}\label{z01}
	&& \E \left[\int_{0}^{T}\left| z_{s}^{n}-z_{s}^{m}\right|
	^{2}ds\right]\\\nonumber
	&&\quad\quad\qquad\leq c\E\left[\int_{0}^{T}|y_s^n-y_s^m|^{\frac{2}{2-\alpha}}ds\right]^{\frac{2-\alpha}{2}}\times\E\left[\int_{0}^{T}\big|f(s,y_{s}^{n},z_{s}^{n})-
	f(s,y_{s}^{m},z_{s}^{m})\big|^{\frac{2}{\alpha}}ds\right]^{\frac{\alpha}{2}}\\\nonumber
	&& \quad\quad\qquad+c \E\left[\sup\limits_{0\leq s\leq T}|\big( L_{s}-y_{s}^{m}\big)^+|^{\frac{p}{p-1}}\right] \E\left[(k^n)^p\right]+c\E\left[\sup\limits_{0\leq s\leq T}|\big( L_{s}-y_{s}^{n}\big)^+|^{\frac{p}{p-1}}\right] \E\left[(k^m)^p\right].
\end{eqnarray}
Then, from Lemma \ref{estf}, \eqref{fat},  \eqref{domconv} and \eqref{leem} (for $\lambda$ large enough and $1<\alpha<2-\frac{2}{e^{\lambda T+1}}$ )
\begin{eqnarray}
	&&\E \left[\int_{0}^{T}\left| z_{s}^{n}-z_{s}^{m}\right|
	^{2}ds\right]\rightarrow 0 \text{ as } (n,m)\rightarrow +\infty.
\end{eqnarray}
Consequently, from \eqref{conse},
$$\E\left[\sup\limits_{0\leq t\leq T}|k^n_t-k^m_t|\right]\rightarrow 0 \text{ as } (n,m)\rightarrow +\infty.$$
Therefore, there exists a pair $(z,k)$ of progressively measurable processes such that
$$\E \left[\int_{0}^{T}\left| z_{s}^{n}-z_{s}\right|
^{2}ds+\sup\limits_{0\leq t\leq T}|k^n_t-k_t|\right]\rightarrow 0 \text{ as } n\rightarrow +\infty.$$
It remains to show that $$\int_{0}^{T}(y_s-L_s)dk_s=0.$$
Clearly, $(k_t)_{0\leq t\leq T}$ is increasing. Moreover, $(y^n,k^n)$ tends to $(y,k)$ uniformly in t in probability. Then
$$\int_{0}^{T}(y^n_s-L_s)dk^n_s\rightarrow\int_{0}^{T}(y_s-L_s)dk_s$$
in probability as $n\rightarrow+\infty$. Therefore, since
$\int_{0}^{T}(y^n_s-L_s)dk^n_s\leq 0,\;n\in\N,$
we have
$ \int_{0}^{T}(y_s-L_s)dk_s\leq 0.$
On the other hand
$ \int_{0}^{T}(y_s-L_s)dk_s\geq 0.$
Thus,
$$ \int_{0}^{T}(y_s-L_s)dk_s= 0.~~\text{a.s.}$$
Hence, $(y,z,k)$ solves the reflected BSDE associated with $(\xi, f, L)$.\qed

We now focus on the uniqueness of the solution for BSDEs \eqref{equ}, Actually, we have the following proposition.
\begin{proposition}
Assume that (\textbf{H.1})-(\textbf{H.4}) are satisfied, then the reflected BSDE associated with $(\xi,f,L,U)$ has at most one solution.
\end{proposition}
\noindent
\textbf{Proof.}
Let suppose that there exists two solutions $(Y,Z,K^+,K^-)$ and $(Y',Z',K'^+,K'^-)$ for \eqref{equ}, and  for $N\in \N^\star$ we set, $ \Delta_{t}:=\left|Y_{t}-Y_{t}'\right| ^{2}+ (A_N)^{-1} $.\\
Following the same argument as in step 2 in the proof of Proposition \ref{penloga}, one can prove that for every $R\in\N$ and for every $\varepsilon>0$  there exists $N_0$ such that for every $N>N_0$
\begin{equation}\label{eps1}
\E \left[\sup_{(T'-\delta')^{+}\leq t \leq T'}\vert
Y_{t}-Y_{t}'\vert^\beta\right]+ \E\left[
\int_{(T'-\delta')^{+}}^{T'}{\left|
Z_{s}-Z_{s}'\right|^{2}\over\left(\vert
Y_{s}'-Y_{s}'\vert^{2}+\nu_{R}\right)^{{2-\beta\over
2}}}ds\right]
\leq \ell e^{C_N\delta'}\E\left[ \vert
Y_{T'}-Y_{T'}'\vert^\beta\right]+\varepsilon.
\end{equation}
where $\nu_{R} = \sup \left\{(A_N)^{-1}, N\geq R
\right\}$ and $\ell$ a universal constant.\\
Taking successively $T'=T$,
$T'=(T-\delta')^+$, $T'=(T-2\delta')^{+}...$ in \eqref{eps1}, we obtain
$$Y=Y',~~Z=Z',~~K^+-K^-=K'^+-K'^-.$$
Finally, let us show that $K^+=K'^+$ and $K^-=K'^-$. For any $t\leq T$, $$\int_{0}^{t}(Y_s-L_s)dK_s=\int_{0}^{t}(Y_s-L_s)dK'_s,$$ where $K=K^+-K^-$ and $K'=K'^+-K'^-.$
But $$\int_{0}^{t}(Y_s-L_s)dK_s=-\int_{0}^{t}(U_s-L_s)dK^-_s \mbox{ and } \int_{0}^{t}(Y_s-L_s)dK'_s=-\int_{0}^{t}(U_s-L_s)dK'^-_s.$$ Then $$ \int_{0}^{t}(U_s-L_s)dK^-_s=\int_{0}^{t}(U_s-L_s)dK'^-_s,~\forall t\leq T.$$ Since $K^-_0=K'^-_0=0$ and $L_t<U_t,~\forall t\leq T$ it follows that $K^-=K'^-$, and we also obtain that $K^+=K'^+$, which completes the proof.\qed \\

After overcoming the main difficulty of this section (Theorem \ref{penloga}) we can now address the question of the existence of a local solution for \eqref{equ}. Actually, we have the following theorem.
\begin{theorem}\label{exl}
	There exists a unique continuous process $Y=(Y_t)_{t\in[0,T]}$ such that
	\begin{itemize}
		\item[(i)] $\E\left[\sup\limits_{s\leq T}|Y_s|^{e^{\lambda T}+1}\right]<+\infty$ and satisfies $L\leq Y\leq U$ and $Y_T=\xi$.
		\item[(ii)] For any stopping time $\tau$ there exists another stopping time $\lambda_{\tau}\geq\tau$, $P$-a.s., and a triplet of processes $(Z^{\tau},K^{\tau,+},K^{\tau,-})\in \mathcal{M}^2\times \mathcal{A}^p\times\mathcal{A}^p,~(K^{\tau,\pm}_{\tau}=0)$ such that $P$-a.s.
		\begin{equation}\label{local}
			\left\{\begin{array}{l}
				{Y}_t ={Y}_{\lambda_{\tau}} + \int_t^{\lambda_{\tau}}
				f(s,{Y}_s,{Z}^{\tau}_s) ds +({K}^{\tau,+}_{\lambda_{\tau}}-{K}^{\tau,+}_t)-({K}^{\tau,-}_{\lambda_{\tau}}-{K}^{\tau,-}_t)\\\qquad- \int_t^{\lambda_{\tau}} {Z}^{\tau}_s dB_s,~~
				t\in[\tau,\lambda_{\tau}];\\
				\int_{\tau}^{\lambda_{\tau}}\left(Y_s-L_s\right)d{K}^{\tau,+}_s=\int_{\tau}^{\lambda_{\tau}}\left(U_s-{Y}_s\right)d{K}^{\tau,-}_s=0.
			\end{array}\right.
		\end{equation}
		\item[(iii)] If $\nu_{\tau}$ and $\pi_{\tau}$ are two stopping times such that:  $$\nu_{\tau}=\inf\{s\geq \tau, Y_s=U_s  \}\wedge T~~\text{and}~~ \pi_{\tau}=\inf\{s\geq \tau, Y_s=L_s  \}\wedge T,$$ then $P$-a.s., $\nu_\tau\vee\pi_\tau\leq\lambda_{\tau}$.
	\end{itemize}
\end{theorem}
\noindent
\textbf{Proof.} After we have proved Theorem \ref{penloga}, the remaining steps to prove Theorem \ref{exl} are actually the same as in  \cite{BHY}. Thus, to avoid repetition, we only give sketch of the poof and for more details we refer the reader to \cite{BHY} pages 914-924.

First, we analyze the following increasing  penalization scheme:  for any $n\geq 0$
\begin{equation}\label{ips}
	\left\{\begin{array}{l}
		Y^n\in \mathcal{S}^{e^{\lambda T}+1},~~Z^n\in\mathcal{M}^2,~~K^{n,-}\in\mathcal{A}^p;\\
		Y^n_t =\xi + \int_t^T \left(
		f(s,Y^n_s,Z^n_s)+n(L_s-Y^n_s)^+\right) ds -(K^{n,-}_T-K^{n,-}_t)- \int_t^T Z^n_s dB_s,\qquad
		t\in[0,T];\\
		Y^n_t\leq U_t,\qquad \forall t\in[0,T]~~~\text{and}~~~\int_{0}^{T}\left(U_s-Y^n_s\right)dK^{n,-}_s=0.
	\end{array}\right.
\end{equation}
First, note that $(Y^n,Z^n,K^{n,-})$ exists due to Theorem \ref{penloga}  and the fact that  $(Y,Z,K)$ is a solution of the reflected BSDE with a lower obstacle associated with $(\xi,f,L)$ iff $(-Y,-Z,K)$ is a solution of the reflected BSDE with an upper obstacle associated with $(-\xi,-f(t,-Y,-Z),-L)$.

Next, since the sequence $f_n(t,y,z)=f(t,y,z)+n(L_t-y)^+ $ is increasing, then from Remark \ref{comprem} we have that for any $n\geq 0$, $Y^n\leq Y^{n+1}\leq U$. Then, $(Y^n_t)_{n\geq 0}$ converges to a lower semi-continuous optional process $Y=(Y_t)_{0\leq t\leq T}$ that satisfies,  $Y_t\leq U_t$, $\forall t\leq T$ $P$-a.s., and $\E\left[\sup\limits_{t\leq T}|Y_t|^{e^{\lambda T}+1}\right]<+\infty.$\\ Next, we put $$\theta_{\tau}^n=\inf\{s\geq \tau, Y^n_s=U_s \}\wedge T,~  \theta_\tau=\lim\limits_{n\rightarrow+\infty}\theta^n_\tau \text{ and } g^n_s= f(s,Y^n_s,Z^n_s)$$ and we show that $Y$ is RCLL on $[\tau,\theta_\tau]$. Indeed, since $K^{n,-}$ does not increase before $\theta_\tau$, $(Y^n_t,Z^n_t)$ satisfy \eqref{ips} with $K^{n,-}_T=K^{n,-}_t=0$ on $[\tau,\theta_\tau]$. Then, as a result of Lemma \ref{estf} and \eqref{fat}, there exists subsequences of $((g^n_s\ind_{[\tau,\theta_{\tau}]}(s))_{s\leq T})_{n\geq 0}$ and $((Z^n_s\ind_{[\tau,\theta_{\tau}]}(s))_{s\leq T})_{n\geq 0}$, which we still index by $n$, and processes $(g_s\ind_{[\tau,\theta_{\tau}]}(s))_{s\leq T}$ and $(Z_s\ind_{[\tau,\theta_{\tau}]}(s))_{s\leq T}$ such that for any stopping time $\bar{\gamma}$ satisfying $\tau\leq \bar{\gamma}\leq \theta_{\tau}$the following weak convergence holds:
$$\int_{\tau}^{\bar{\gamma}}Z^n_sdB_s\rightharpoonup \int_{\tau}^{\bar{\gamma}}Z_sdB_s~~~\text{and}~~~ \int_{\tau}^{\bar{\gamma}}g^n_sds \rightharpoonup\int_{\tau}^{\bar{\gamma}}g_sds, ~~~\text{as }~ n\rightarrow+\infty.$$ It follows that $$K^{n,+}_{{\bar{\gamma}}}\rightarrow K^+_{{\bar{\gamma}}} \text{ and } Y_t=Y_{\tau}-\int_{\tau}^{t}g_sds-K^+_t+\int_{\tau}^{t}Z_s dB_s,$$
 such that $\E\left[(K^+_{\theta_\tau})^p\right]\leq \underset{n\rightarrow+\infty}{\lim\inf}~\E\left[(K^{n,+}_{\theta_{\tau}})^p\right]<+\infty$. Since $Y^n\leq Y^{n+1}$, we can deduce from a result by S. Peng [\cite{PE}, Lemma 2.2] that $Y$ is RCLL on $[\tau,\theta_\tau]$. Next, we can show as in \cite{BHY} that we have the following proposition which can be considered as steps of the proof.
 \begin{proposition}\label{pro} Assume that (\textbf{H.1})-(\textbf{H.4}) are satisfied. Then, the following holds true:
 	\begin{itemize}
 		\item[(i)] $ P$-a.s.,  $Y_{\theta_{\tau}}\ind_{\{\theta_{\tau}< T\}}=U_{\theta_{\tau}}\ind_{\{\theta_{\tau}< T\}}$ and $ P$-a.s.,  $\forall t\leq T$, $L_t\leq Y_t.$
 		\item[(ii)]
 		There exist two adapted processes $(\bar{K}^{\tau,+}_t)_{0\leq t\leq T}$ and $(\bar{Z}^{\tau}_t)_{0\leq t\leq T}$ such that  $(Y_t,\bar{Z}^{\tau}_t,\bar{K}^{\tau,+}_t,0)_{0\leq t\leq T}$ is a local solution of the reflected BSDE \eqref{equ} on $[\tau,\theta_{\tau}]$, which mean it satisfies the following:
 		\begin{equation}\label{loc}
 			\left\{\begin{array}{l}
 				\bar{Z}^{\tau}\in\mathcal{M}^2,~~\bar{K}^{\tau,+}\in\mathcal{A}^p;\\
 				Y_t =Y_{\theta_{\tau}} + \int_t^{\theta_{\tau}}
 				f(s,Y_s,\bar{Z}^{\tau}_s) ds +(\bar{K}^{\tau,+}_{\theta_{\tau}}-\bar{K}^{\tau,+}_t)- \int_t^{\theta_{\tau}} \bar{Z}_s dB_s,\qquad
 				\forall t\in[\tau,\theta_{\tau}];\\
 				Y_T=\xi;\\
 				\forall t\in[\tau,\theta_{\tau}]\qquad L_t\leq Y_t\leq U_t,~~~\text{and}~~~\int_{\tau}^{\theta_{\tau}}\left(Y_s-L_s\right)d\bar{K}^{\tau,+}_s=0.
 			\end{array}\right.
 		\end{equation}
 		\item[(iii)] $v_{\tau}=\inf\{s\geq\tau,Y_s=U_s\}\wedge T$, then $v_{\tau}\leq\theta_{\tau}$.
 	\end{itemize}
 \end{proposition}
Now by analyzing the decreasing penalization scheme that is: for any $m\geq 0$
\begin{equation}\label{dps}
	\left\{\begin{array}{l}
		\E\left[\sup\limits_{0\leq s\leq T}|\tilde{Y}^m_s|^{e^{\lambda T}+1}+\int_{0}^{T}|\tilde{Z}^m_s|^2 ds+({K}^{m,+}_{T})^p\right]<+\infty;\\
		\tilde{Y}^m_t =\xi + \int_t^T \left(
		f(s,\tilde{Y}^m_s,\tilde{Z}^m_s)-m(\tilde{Y}^m_s-U_s)^+\right) ds +(K^{m,+}_T-K^{m,+}_t)- \int_t^T \tilde{Z}^m_s dB_s,\qquad
		t\in[0,T];\\
		\tilde{Y}^m_t\geq L_t,\qquad t\in[0,T]~~~\text{and}~~~\int_{0}^{T}\left(\tilde{Y}^m_s-L_s\right)dK^{m,+}_s=0,
	\end{array}\right.
\end{equation}
($(\tilde{Y}^m,\tilde{Z}^m,K^{m,+})$ exists due to Theorem \ref{penloga}) we can also show that
\begin{proposition}\label{pro2}
	The following hold:
	\begin{itemize}
		\item[(i)] $ P$-a.s.,
		$\tilde{Y}_{\delta_{\tau}}\ind_{\{\delta_{\tau}< T\}}=L_{\delta_{\tau}}\ind_{\{\delta_{\tau}< T\}}$ and $P$-a.s., $\forall t\leq T$, we have $\tilde{Y}_t\leq U_t$.
		\item[(ii)] There exist a pair of adapted processes $(\tilde{Z}^{\tau}_t,\tilde{K}^{\tau,-}_t)_{t\leq T}$ such that the quadruple $(\tilde{Y}_t,\tilde{Z}^{\tau}_t,0,\tilde{K}^{\tau,-}_t)_{t\leq T}$ satisfy
		\begin{equation}\label{loc1}
			\left\{\begin{array}{l}
				\tilde{Z}^{\tau}\in\mathcal{M}^2,~~\tilde{K}^{\tau,-}\in\mathcal{A}^p;\\
				\tilde{Y}_t =\tilde{Y}_{\delta_{\tau}} + \int_t^{\delta_{\tau}}
				f(s,\tilde{Y}_s,\tilde{Z}^{\tau}_s) ds -(\tilde{K}^{\tau,-}_{\delta_{\tau}}-\tilde{K}^{\tau,-}_t)- \int_t^{\delta_{\tau}} \tilde{Z}^{\tau}_s dB_s,\qquad
				t\in[\tau,\delta_{\tau}];\\
				\tilde{Y}_T=\xi;\\
				L_t\leq \tilde{Y}_t\leq U_t,\quad \forall t\in[\tau,\delta_{\tau}]~\text{ and }~\int_{\tau}^{\delta_{\tau}}\left(U_s-\tilde{Y}_s\right)d\tilde{K}^{\tau,-}_s=0.
			\end{array}\right.
		\end{equation}
		\item[(iii)] Put $\mu_{\tau}=\inf\{s\geq\tau,\tilde{Y}_s=L_s\}\wedge T$, then $\mu_{\tau}\leq \delta_{\tau}$,
	\end{itemize}
where $\tilde{Y}=\lim\limits_{m\rightarrow+\infty} \tilde{Y}^m$ and $\delta_{\tau}=\lim\limits_{m\rightarrow+\infty}\delta^m_{\tau}$ with $\delta_{\tau}^m=\inf\{s\geq \tau, \tilde{Y}^m_s=L_s \}\wedge T,~~\forall m\geq 0$.
\end{proposition}
Next using the comparison result and the technique in \cite{BHY} page 923, we can prove that $P$-a.s., for any $t\leq T$, $Y_t=\tilde{Y}_t$. Finally, we proceed, once again, as in \cite{BHY} page 924 to finish the proof.\\

Next we can proceed as in [\cite{hahs}, Theorem 3.7] to show that the local solution is actually a global one, meaning, we have the  following theorem:
\begin{theorem}\label{eu}
 Under the assumptions (\textbf{H.1}), (\textbf{H.2}), (\textbf{H.3}) and (\textbf{H.4}),
 the reflected BSDE \eqref{equ} associated with $(\xi,f,L,U)$ has a unique solution that is the quadruple $(Y,Z,K^+,K^-)$.
\end{theorem}
\section{Mixed zero-sum stochastic differential game problem}
Now we deal with an application of the double barrier reflected BSDEs tool for solving stochastic mixed games problems. First, let us briefly describe the setting of the considered problem. In the sequel $\Omega={\cal C}([0,T],\R^d)$ is the space of continuous functions from $[0,T]$ to $\R^d$.

Put $||\omega||_t=\sup\limits_{s\leq t}|\omega_s|$ and let us consider a mapping $\sigma:\, (t,\omega)\in[0,T]\times \Omega
\rightarrow \sigma(t,\omega)\in\R^{d}\bigotimes \R^{d}$ satisfying the
following assumptions
\begin{itemize}
	\item[\textbf{(A1)}]
\begin{itemize}
\item[\textbf{(i)}]$\sigma$ is {$\mathcal{P}$}-measurable and invertible.
\item[\textbf{(ii)}]  There exists a constant $C>0$ such that: $\forall (t,\omega,\omega')\in[0,T]\times\Omega\times\Omega$
$$|\sigma(t,\omega)-\sigma(t,\omega')|\leq C||\omega-\omega'||_t,~ |\sigma(t,\omega)|\leq
C(1+||\omega||_t) \text{ and } |\sigma^{-1}(t,\omega)|\leq
C.$$
\end{itemize}
\end{itemize}
\medskip\noindent
Let $x_0\in\R^d$ and $x=(x_t)_{t\leq T}$ be the solution of the
following standard functional differential equation:
\begin{equation}\label{EDS}
x_t=x_0+\int_0^t\sigma(s,x)dB_s,\quad t\leq T;
\end{equation}
The assumptions on $\sigma$ imply that the equation \eqref{EDS} has a unique solution $x$ (see, \cite{RY} page 375). Moreover,
 \begin{equation}\label{estim-eds}
 \E[(||x||_T)^n]<+\infty,\quad \forall n\in[1,+\infty[~~(\cite{KS},\text{ pp.
 306}).
 \end{equation}
Let us consider now a compact metric space $\bar{U}$ (resp. $V$) and $\mathcal{U}$ (resp. $\mathcal{V}$) the space
of all $\mathcal{P}$-measurable processes with values in $\bar{U}$ (resp. $V$), and let
$\varphi:[0,T]\times \R^d \times \bar{U}\times V\rightarrow \R^d$ and $h:[0,T]\times \R^d \times \bar{U}\times V\rightarrow \R^d$ be such
that:
\begin{itemize}
	\item[\textbf{(A2)}]
\begin{itemize}
\item[\textbf{(i)}] For each $(u,v)\in \bar{U}\times V$, the function $(t,x)\rightarrow \varphi(t,x,u,v)$
is predictable.
\item[\textbf{(ii)}] $\forall(t,x)\in[0,T]\times\R^d$, $\varphi(t,x,.,.)$ and $h(t,x,.,.)$ are continuous on $\bar{U}\times V$.
\item[\textbf{(iii)}] There exists a real constant $K>0$ such that
\begin{equation}\label{generateur-cond}
|h(t,x,u,v)|+|\varphi(t,x,u,v)|\leq K(1+||x||_t),\quad \forall (t,x,u,v)\in [0,T]\times\R^d\times \bar{U}\times V.
\end{equation}
\end{itemize}
\end{itemize}
Under the previous assumption, and for any $(u,v)\in\mathcal{U}\times\mathcal{V}$, we define a probability on ($\Omega,\mathcal{F}$) by

$$\frac{dP^{(u,v)}}{dP}= \exp\bigg\{\int_0^T\sigma^{-1}(s,x)\varphi(s,x,u_s,v_s)dB_s-
\frac{1}{2}\int_0^T|\sigma^{-1}(s,x)\varphi(s,x,u_s,v_s)|^2ds\bigg\}.$$
\par\noindent
We now consider the payoff
\begin{equation}\label{pay}
J(u,\tau;v,\sigma)=\E^{(u,v)}\left[\int_{0}^{\tau\wedge\sigma}h(s,x,u_s,v_s)ds+L_\sigma\ind_{\{ \sigma\leq\tau<T \}}+U_\tau\ind_{\{ \tau<\sigma\}}+\xi\ind_{\{ \tau\wedge\sigma=T \}}\right],
\end{equation}
where $L$, $U$ and $\xi$ are those of the previous sections.
The problem we are interested in is finding a saddle-point for the payoff functional $J(u,\tau;v,\sigma)$, meaning we are looking for two interventions strategies $(u^*,\tau^*)$ and $(v^*,\sigma^*)$ that satisfies
\begin{equation}\label{j}
J(u^*,\tau^*;v,\sigma)\leq J(u^*,\tau^*;v^*,\sigma^*)\leq J(u,\tau;v^*,\sigma^*).
\end{equation}
Now we define the Hamiltonian associated with this mixed stochastic game problem by
$$ H(t,x,z,u,v):= z\sigma^{-1}(t,x)\varphi(t,x,u,v)+h(t,x,u,v)\quad \forall(t,x,z,u,v)\in[0,T]\times\R^{d}\times\R^{d}\times \bar{U}\times V.$$
Under Isaacs's condition, and Benes' theorem \cite{B}, there exits a couple of $\mathcal{P}\otimes\mathcal{B}$-measurable functions $u^*\equiv u^*(t,x,z)$ and $v^*\equiv v^*(t,x,z)$ with values in $\bar{U}$ and $V$ respectively such that: $\forall (t,x,u,v)\in[0,T]\times\R ^d\times \bar{U}\times V$
$$H^*(t,x,z)=H(t,x,z,u^*(t,x,z),v^*(t,x,z))=\inf_{u\in \bar{U}}\sup_{v\in V}H(t,x,z,u,v)=\sup_{v\in V}\inf_{u\in \bar{U}}H(t,x,z,u,v).$$
\begin{theorem}\label{the5.1}
Under the assumptions (\textbf{A1}) and (\textbf{A2}), there exists a quadruple of adapted processes $(Y^*,Z^*,K^{*,+},K^{*,-})$ such that it is the uniques solution of the finite horizon reflected BSDE associated with $(\xi,H^*,L,U)$. We denote by $\tau^*$ and $\sigma^*$ the stopping times defined as follow:
$$ \sigma^*=\inf\{t\geq 0,~Y^*_t =  L_t\}\wedge T ~~\text{and}~~ \tau^*=\inf\{t\geq 0,~Y^*_t = U_t\}\wedge T.$$
Then, $Y^*_0=J(u^*,\tau^*;v^*,\sigma^*)$ and  $(u^*,\tau^*;v^*,\sigma^*)$
is a saddle-point for the mixed stochastic game problem.
\end{theorem}
\noindent
\textbf{Proof.} Since $H^*$ satisfy \textbf{(H.3)} and \textbf{(H.4)} (see \cite{EO}),  the quadruple $(Y^*,Z^*,K^{*,+},K^{*,-})$ exists and unique. the rest of the proof is classical, thus we leave to the reader.

\section{Connection with Double Obstacle Variational Inequalities}
Let $b:[0,T]\times \R^d\rightarrow \R^d , \sigma:[0,T]\times \R^d\rightarrow \R^{d\times d} $ be two globally Lipschitz functions and let us consider the following SDE:
$$dX_t=b(t,X_t)dt+\sigma(t,X_t)dB_t,~~t\leq T.$$
We denote by $ (X^{t,x}_s)_{s \geq t}$ the unique solution of the previous SDE starting from $x$ at time $s = t$. \\
Now we are given four functions
$$f:[0,T]\times \R^d\times\R\times\R^d\rightarrow \R,~~g:\R^d\rightarrow\R, ~\text{and}~h,h':[0,T]\times\R^d\rightarrow\R,$$
such that the followings holds:
\begin{itemize}
\item[(\textbf{H'.1})] $f$ satisfies assumptions \textbf{(H.3)} and \textbf{(H.4)}, moreover there exists $p>1$ such that for every $(t,x)\in[0,T]\times\R^d$
$$\E\left[\int_{0}^{T}|f(s,X^{t,x}_s,0,0)|^pds\right]<+\infty,$$
\item [(\textbf{H'.2})] $\forall (t,x)\in[0,T]\times\R^d$, $h(t,x)<h'(t,x)$ and $h(T,x)\leq g(x)\leq h'(T,x)$, in addition there exists a constant $C>0$ such that
$$|h'(t,x)|+|h(t,x)|+|g(x)|\leq C(1+||x||_t).$$
\end{itemize}
\subsection{Connection with one Obstacle variational inequalities}

Let us define $(Y^{t,x}_s,Z^{t,x}_s,K^{t,x}_s)_{s\in[t,T]}$  the solution of the following reflected BSDE :
\begin{equation}
	\left\{
	\begin{array}{l}
		Y^{t,x}_s = g(X^{t,x}_T) + \int_s^T
		f(u,X^{t,x}_u,Y^{t,x}_u,Z^{t,x}_u) du +K^{t,x}_T-K^{t,x}_s- \int_s^T Z^{t,x}_u dB_u\\
		\forall s\in[t,T]~~~~ Y^{t,x}_s\geq h(s,X^{t,x}_s),~~\text{and}~~ \int_{t}^{T} \left(h(u,X^{t,x}_u)- Y^{t,x}_u\right)dK^{t,x}_u=0
	\end{array}\right.
\end{equation}
Moreover on $[0,t]$ we set $Y^{t,x}_s=Y^{t,x}_t$, $Z^{t,x}_s=K^{t,x}_s=0.$\\
For every $(t,x)$ we will show that $Y^{t,x}_t$ is  deterministic and we define a function
\begin{equation}
	u(t,x)=Y^{t,x}_t.
\end{equation}
First we will prove that $u$ is continuous and is a viscosity solution of the following obstacle problem. $\forall(t,x)\in[0,T]\times \R ^d$
\begin{equation}\label{kkk7}
	\min\left[u(t,x)-h(t,x),\\
	-\frac{\partial u}{\partial t}(t,x)-\mathcal{L}u(t,x)-f(t,x,u(t,x),\sigma(t,x)\nabla u(t,x))\right]=0,
\end{equation}
with $ u(T,x)=g(x),\; x\in\R^d.$
Then we will prove that it is the unique continuous viscosity solution that belongs to some class of functions.
\subsubsection{Continuity}
\begin{theorem}\label{conti7}
	For every $(t,x)\in[0,T]\times\R^d$, the function $ u(t,x)=Y^{t,x}_t$ is continuous and of polynomial growth.
\end{theorem}

\noindent
\textbf{Proof.} Let $(t_n,x_n)\rightarrow(t,x)$, since $|Y^{t,x}_t-Y^{t_n,x_n}_{t_n}|$ is deterministic, we have
\begin{eqnarray*}
	&& |Y^{t,x}_t-Y^{t_n,x_n}_{t_n}|=\E \left(|Y^{t,x}_t-Y^{t_n,x_n}_{t_n}|\right)\\ \nonumber
	&&\qquad\qquad\qquad~~\leq \E \left(|Y^{t,x}_t-Y^{t,x}_{t_n}|\right)+\E \left(|Y^{t,x}_{t_n}-Y^{t_n,x_n}_{t_n}|\right).
\end{eqnarray*}
Then from \eqref{estKYZ}  and Lemma \ref{estf} we get $\lim\limits_{n\rightarrow +\infty} \E \left[|Y^{t,x}_t-Y^{t,x}_{t_n}|\right]=0.$\\

Now we shall show that $\lim\limits_{n\rightarrow +\infty} \E \left[|Y^{t,x}_{t_n}-Y^{t_n,x_n}_{t_n}|\right]=0,$ for that we use the fact that $$\E \left[|Y^{t,x}_{t_n}-Y^{t_n,x_n}_{t_n}|\right]\leq \E \left[\sup_{0\leq s\leq T}|Y^{t,x}_{s}-Y^{t_n,x_n}_{s}|\right].$$ We proceed as in the proof of step 2 Proposition \ref{penloga}, which that for $\beta\in]1,\min\left(3-\alpha,2\right)[,$ we have:
\begin{align*}
	& \E \left[\sup_{(T'-\delta')^{+}\leq s \leq T'}|Y^{t,x}_{s}-Y^{t_n,x_n}_{s}|^\beta\right] + \E\left[
	\int_{(T'-\delta')^{+}}^{T'}\dfrac{\left|
		Z_{s}^{t,x}-Z_{s}^{t_n,x_n}\right|^{2}}{\left(\vert
		Y_{s}^{t,x}-Y_{s}^{t_n,x_n}\vert^{2}+ \nu_{R}
		\right)^{\frac{2-\beta}{2}}}\right]ds
	\\ & \qquad\qquad\leq
	\ell e^{C_N\delta'}\E\left[ \vert
	g\left(X^{t,x}_{T'}\right)-g\left(X^{t_n,x_n}_{T'}\right)\vert^\beta\right]+ \ell
	\left(\dfrac{A_N^{\frac{2M^2\delta'\beta}{\beta-1}}}{(A_N)^{\frac{\beta}{2}}}+\dfrac{A_N^{\frac{2M^2\delta'\beta}{\beta-1}}}{(A_N)^{\frac{\kappa}{r}}}\right)
	\\ &\qquad\qquad +\ell e^{C_N\delta'}\beta [2N^2+\nu_1
	]^{\frac{\beta-1}{2}}
	\\&\qquad\qquad\times \E \left[\int_{t}^{T'}|f(u,X^{t,x}_u,Y^{t,x}_u,Z^{t,x}_u)-f(u,X^{t_n,x_n}_u,Y^{t,x}_u,Z^{t,x}_u)|du\right]
	\\&\qquad\qquad+ \ell e^{C_N\delta'} \left(\E\left[\sup_{0\leq s\leq T}|\left(|h(s,X^{t,x}_s)-h(s,X^{t_n,x_n}_s)|^2+(A_N)^{-1}\right)^{\frac{\beta}{2}-1}\right.\right.
	\\&\qquad\qquad\times\left.\left.{\left(h(s,X^{t,x}_s)-h(s,X^{t_n,x_n}_s)\right)|^{\frac{p}{p-1}}} \right] \right)^{\frac{p-1}{p}}.
\end{align*}
Since $f$, $g$ and $h$ are continuous in the $x$-variable, then for $ \delta' < (\beta
-1)\min\left(\frac{1}{4M^2},\frac{\kappa}{2rM^2\beta}\right)$ we pass to the limit on $n$ and then to the limit on $N$, and by taking successively $T'=T$,
$T'=(T-\delta')^+$, $T'=(T-2\delta')^{+}...$ we get for every $\beta\in]1,\min\left(3-\alpha,2\right)[,$
$$\lim_{n\rightarrow +\infty}\E \left[\sup_{0\leq s\leq T}\vert
Y_{s}^{t,x}-Y_{s}^{t_n,x_n}\vert^\beta\right]=0.$$
Finally, since $\beta>1$, then by using H\"older's inequality the result follows.
The polynomial growth of $u$ follows from \eqref{estKYZ}.\qed
\subsubsection{Existence of the solution}
\begin{theorem}
	Assume that (\textbf{H'.1}) and (\textbf{H'.2}) are satisfied, then the function $u:(t,x)\rightarrow u(t,x)=Y^{t,x}_t$ is a viscosity solution
	of the obstacle problem in (\ref{kkk7}).
\end{theorem}

\noindent
\textbf{Proof.}
Let us consider the following reflected BSDE:
\begin{equation}\label{pen}
	Y^{t,x,n}_s = g(X^{t,x}_T) + \int_s^T
	f_n(r,X^{t,x}_r,Y^{t,x,n}_r,Z^{t,x,n}_r) dr - \int_s^T Z^{t,x,n}_r dB_r,
\end{equation}
where $$f_n(r,X^{t,x}_r,Y^{t,x,n}_r,Z^{t,x,n}_r)=f(r,X^{t,x}_r,Y^{t,x,n}_r,Z^{t,x,n}_r)+n\left(Y^{t,x,n}_r-h(r,X^{t,x}_r)\right)^-.$$
Then, from \cite{BKK},  $u_n(t,x)=Y^{t,x,n}_t$ is the viscosity solution of
\begin{equation}\label{viscon}
	\frac{\partial u_n}{\partial t}(t,x)+\mathcal{L}u_n(t,x)+f_n(t,x,u_n(t,x),\sigma(t,x)\nabla u_n(t,x))=0.
\end{equation}
From the comparison theorem  we have that $u_n$ is increasing, and we can argue as in \cite{EKPPQ} to show that $u_n$ converges to $u$ the solution of \eqref{kkk7}.\qed
\subsection{Connection with Double Obstacle Variational Inequalities}
Let  $(Y^{t,x}_s,Z^{t,x}_s,K^{+,t,x}_s,K^{-,t,x}_s)_{t\leq s\leq T}$ as the solution of the following reflected BSDE :
\begin{equation}
\left\{
\begin{array}{l}
Y^{t,x}_s = g(X^{t,x}_T) + \int_s^T
f(u,X^{t,x}_u,Y^{t,x}_u,Z^{t,x}_u) du+\int_{s}^{T}dK^{+,t,x}_u-\int_{s}^{T}dK^{-,t,x}_u- \int_s^T Z^{t,x}_u dB_u,\\
\forall s\in[t,T],~~~~h(s,X^{t,x}_s)\leq Y^{t,x}_s\leq h'(s,X^{t,x}_s),\\
\int_{t}^{T} \left( Y^{t,x}_u-h(u,X^{t,x}_u)\right)dK^{+,t,x}_u=\int_{t}^{T} \left(h'(u,X^{t,x}_u)- Y^{t,x}_u\right)dK^{-,t,x}_u=0.
\end{array}\right.
\end{equation}
The objective of this section is to show that $u(t,x)=Y^{t,x}_t$ is continuous and it is the solution in the viscosity sense of the following obstacle problem:
\begin{equation}\label{kkk}
	\left\{
	\begin{array}{l}
\min\left[u(t,x)-h(t,x),\max\bigg\{
-\frac{\partial u}{\partial t}(t,x)-\mathcal{L}u(t,x)\right.\\
	\qquad\qquad\left.-f(t,x,u(t,x),\sigma(t,x)\nabla u(t,x)),u(t,x)-h'(t,x)\bigg\}\right]=0;~~~(t,x)\in[0,T)\times \R ^d\\
u(T,x)=g(x),~~\forall x\in\R ^d.	
	\end{array}\right.
\end{equation}
\subsection{The continuity of the viscosity solution}
\begin{proposition}\label{conti}
For every $(t,x)\in[0,T]\times\R^d$, the function $ u(t,x)=Y^{t,x}_t$ is continuous.
\end{proposition}

\noindent
\textbf{Proof.} For any $n\geq 0$ let $(   \underline{Y}^{t,x,n}_s)_{s\leq T}$ (resp. $(\overline{Y}^{t,x,n}_s)_{s\leq T}$) be the first component of the unique solution of the BDSE with one reflecting lower (resp. upper) barrier associated with $(g(X^{t,x}_T),f(s,X^{t,x}_s,y,z)-n(h'(s,X^{t,x}_s)-y)^-,h(s,X^{t,x}_s))$ (resp. $(g(X^{t,x}_T),f(s,X^{t,x}_s,y,z)+n(h(s,X^{t,x}_s)-y)^+,h'(s,X^{t,x}_s))$). As it has been shown in the previous subsection, for any $n\geq 0$ there exist two deterministic functions $\underline{u}^n(t,x)=\underline{Y}^{t,x,n}_t$ and $\bar{u}^n(t,x)=\overline{Y}^{t,x,n}_t$ such that they are the viscosity solution of
\begin{eqnarray}\label{low}
&&\min\left[u(t,x)-h(t,x),
-\frac{\partial u}{\partial t}(t,x)-\mathcal{L}u(t,x)-f(t,x,u(t,x),\sigma(t,x)\nabla u(t,x))\right. \\\nonumber
&& \qquad \quad +n(h'(t,x)-u(t,x))^-\left.\right]=0,~~(t,x)\in[0,T)\times\R^d;~~u(T,x)=g(x),
\end{eqnarray}
and
\begin{eqnarray}\label{upp}
&&\max\left[u(t,x)-h'(t,x),
-\frac{\partial u}{\partial t}(t,x)-\mathcal{L}u(t,x)-f(t,x,u(t,x),\sigma(t,x)\nabla u(t,x))\right. \\\nonumber
&& \qquad \quad -n(h(t,x)-u(t,x))^+\left.\right]=0,~~(t,x)\in[0,T)\times\R^d;~~u(T,x)=g(x)
\end{eqnarray}
respectively. Now thanks to the results of the previous sections, the sequence $(\underline{Y}^{t,x,n})_{n\geq 0}$ converges increasingly to $Y^{t,x}$ and the sequence $(\overline{Y}^{t,x,n})_{n\geq 0}$ converges decreasingly to the same $Y^{t,x}$, meaning that $\underline{u}^n(t,x)\searrow u(t,x)$ and $\bar{u}^n(t,x)\nearrow u(t,x)$. Since $\underline{u}^n$ and $\bar{u}^n$ are both continuous, then $u$ is at the same time, lower and upper semi-continuous, therefore, it is continuous.
\subsection{Existence of the solution}
\begin{theorem}
Assume that (\textbf{H'.1}) and (\textbf{H'.2}) are satisfied, then the function $u:(t,x)\mapsto u(t,x)=Y^{t,x}_t$ is a viscosity solution
of the obstacle problem in (\ref{kkk}).
\end{theorem}

\noindent
\textbf{Proof.} First note that since $\underline{u}^n$, $\bar{u}^n$ and $u$ are continuous, then from Dini's Lemma,  $\underline{u}^n$ and $\bar{u}^n$ converge uniformly to $u$ on compact subsets of $[0,T]\times\R^d$.\\
Let us now show that $u$ is a viscosity subsolution of \eqref{kkk}. Let  $\phi$ be a $ \mathcal{C}^{1,2}\left((0,T)\times\R^d\right)$, and
$(t_n,x_n)$ be a sequence of local maximum points of $\underline{u}^n-\phi$ such that it converges to $(t,x)$. For $n$ large enough we have $\underline{u}^n(t_n,x_n)>h(t_n,x_n)$, and since $\underline{u}^n$ is a viscosity solution of \eqref{low} we have:
\begin{eqnarray*}
&&-\frac{\partial \phi}{\partial t}(t_n,x_n)-\mathcal{L}\phi(t_n,x_n)
-f(t_n,x_n,\underline{u}^n(t_n,x_n),\sigma(t_n,x_n)\nabla \phi(t_n,x_n))\\&&
\qquad\qquad\qquad+n(h'(t_n,x_n)-\underline{u}^n(t_n,x_n))^-\leq 0,
\end{eqnarray*}
then,
$$-\frac{\partial \phi}{\partial t}(t_n,x_n)-\mathcal{L}\phi(t_n,x_n)-f(t_n,x_n,\underline{u}^n(t_n,x_n),\sigma(t_n,x_n)\nabla \phi(t_n,x_n))\leq 0.$$
Now due to the continuity of the functions and the uniform convergence of $\underline{u}^n$ we obtain
$$-\frac{\partial \phi}{\partial t}(t,x)-\mathcal{L}\phi(t,x)-f(t,x,u(t,x),\sigma(t,x)\nabla \phi(t,x))\leq 0.$$
Since $u(T,x)=g(x)$ and $h(t,x)\leq u(t,x)\leq h'(t,x)$, $u$ is a viscosity subsolution of \eqref{kkk}. In the same way, with converse inequalities, we show that $u$ is also a viscosity supersolution of \eqref{kkk}.
\subsection{Uniqueness of the viscosity solution}
We are going now to address the
question of uniqueness of the viscosity solution of \eqref{kkk}. But first we recall the following proposition.
\begin{proposition}
$w$ is a viscosity solution of
\begin{equation}\label{w}
\left\{\begin{array}{l}
\min\left[w(t,x)-h(t,x),
-\frac{\partial w}{\partial t}(t,x)-\mathcal{L}w(t,x)\right.\\
\left.-f(t,x,w(t,x),\sigma(t,x)\nabla w(t,x))\right]=0,~~(t,x)\in[0,T[\times\R^d\\
w(T,x)=g(x),~~x\in\R^d,
\end{array}\right.
\end{equation}
iff $\overline{w}(t,x) =
e^tw(t,x)$, for any $t\in[0,T]$ and $x\in \mathbb{R}^d$, is a viscosity solution  of
\begin{equation}\label{tra}
\left\{\begin{array}{l}
\min\left[\overline{w}(t,x)-e^th(t,x),-\frac{\partial \overline{w}}{\partial t}(t,x)+\overline{w}(t,x)-\mathcal{L}\overline{w}(t,x)\right.\\
\left.-e^t f(t,x,e^{-t}\overline{w}(t,x),\sigma(t,x)\nabla (e^{-t}\overline{w}(t,x)))\right]=0,~~(t,x)\in[0,T[\times \R^d\\
\overline{w}(T,x)=e^{T}g(x),~~x\in\R^d.
\end{array}\right.
\end{equation}
\end{proposition}

We now have the following theorem.
\begin{theorem}
Under (\textbf{H'.1}) and (\textbf{H'.2}), the equation \eqref{kkk} has at most one solution.
\end{theorem}

\noindent
\textbf{Proof.} In order to prove the uniqueness of the solution it is enough to show that if $v$ and $u$ are viscosity supersolution and subsolution of \eqref{kkk} respectively, then $$u(t,x)\leq v(t,x),~~~\forall (t,x)\in[0,T]\times \R^d.$$
First, note that $v\geq h$ et $u\leq h'$ and set $\bar{v}:=v\wedge h'$ and $\bar{u}:=u\vee h$. Then, $\bar{u}$ (resp. $\bar{v}$) is a viscosity subsolution (resp. supersolution)  of \eqref{kkk}. It follows that $\bar{u}$ (resp. $\bar{v}$) is a viscosity subsolution (resp. supersolution) of \eqref{kkk7} 

Now we show that $\bar{v}$ and $\bar{u}$ satisfy $\bar{u}\leq \bar{v}$. 
Actually for some $R>0$ suppose there exists
$(\overline{t},\overline{x})\in[0,T]\times B_R$
$(B_R := \{x\in \R^d; |x|<R\})$ such that:
\begin{equation}
\label{comp-equ1}
\max\limits_{t,x}(u'(t,x)-v'(t,x))=u'(\overline{t},\overline{x})-v'(\overline{t},\overline{x})=\eta>0,
\end{equation}
where $v'(t,x)=e^t\bar{v}(t,x)$ and $u'(t,x)=e^t\bar{u}(t,x)$ for any $t\in[0,T]$ and $x\in \R^d$.

Let us take
$\theta$, $\lambda$ and $\beta \in (0,1]$ small enough.
 Then, for a small $\epsilon>0$, let us define:
\begin{equation}
\label{phi}
\Phi_{\epsilon}(t,x,y)=(1-\lambda)u'(t,x)-v'(t,y)-\frac{1}{2\epsilon}|x-y|^{4}
-\theta( |x-\overline{x}|^{4}+|y-\overline{x}|^{4})-\beta
(t-\overline{t})^2.
\end{equation}
We have $u'$ and $v'$ are bounded, there exists a
$(t_{\epsilon},x_{\epsilon},y_{\epsilon})\in [0,T]\times B_R
\times B_R $, for $R$ large enough, such that:
$$\Phi_{\epsilon}(t_{\epsilon},x_{\epsilon},y_{\epsilon})=\max\limits_{(t,x,y)}\Phi_{\epsilon}(t,x,y).$$
On the other hand, from
$2\Phi_{\epsilon}(t_{\epsilon},x_{\epsilon},y_{\epsilon})\geq
\Phi_{\epsilon}(t_{\epsilon},x_{\epsilon},x_{\epsilon})+\Phi_{\epsilon}(t_{\epsilon},y_{\epsilon},y_{\epsilon})$,
we have
\begin{equation*}
\frac{1}{\epsilon}|x_{\epsilon} -y_{\epsilon}|^{4} \leq
(1-\lambda)(u'(t_{\epsilon},x_{\epsilon})-u'(t_{\epsilon},y_{\epsilon}))+(v'(t_{\epsilon},x_{\epsilon})-v'(t_{\epsilon},y_{\epsilon})),
\end{equation*}
and consequently $\frac{1}{\epsilon}|x_{\epsilon}
-y_{\epsilon}|^{4}$ is bounded, and as $\epsilon\rightarrow 0$,
$|x_{\epsilon} -y_{\epsilon}|\rightarrow 0$. Since $u'$ and $v'$ are
uniformly continuous on $[0,T]\times \overline{B}_R$, then
$\frac{1}{2\epsilon}|x_{\epsilon} -y_{\epsilon}|^{4}\rightarrow 0$
as
$\epsilon\rightarrow 0.$\\
Since
 $$(1-\lambda)u'(\overline{t},\overline{x})-v'(\overline{t},\overline{x}) \leq
 \Phi_{\epsilon}(t_{\epsilon},x_{\epsilon},y_{\epsilon})\leq (1-\lambda)u'(t_{\epsilon},x_\epsilon)-v'(t_{\epsilon},y_\epsilon),$$
it follow as $\lambda\rightarrow 0$ and the continuity of $u'$ and
$v'$ that, up to a subsequence,
 \begin{equation}\label{subsequence}
 (t_\epsilon,x_\epsilon,y_\epsilon)\rightarrow (\overline{t},\overline{x},\overline{x}).
 \end{equation}
 Next let us show that $t_{\epsilon} <T.$ Actually if $t_{\epsilon}
 =T$ then,
 $$
 \Phi_{\epsilon}(\overline{t},\overline{x},\overline{x})\leq
 \Phi_{\epsilon}(T,x_{\epsilon},y_{\epsilon}),$$ and,
 $$
 (1-\lambda)u'(\overline{t},\overline{x})-v'(\overline{t},\overline{x})\leq
 (1-\lambda)e^{T}g(x_{\epsilon}) -e^{T}g(y_{\epsilon})- \beta
 (T-t_{\epsilon})^2,
 $$
 since $u'(T,x_{\epsilon})=e^{T}g(x_{\epsilon})$,
 $v'(T,y_{\epsilon})=e^{T}g(y_{\epsilon})$ and $g$ is uniformly
 continuous on $\overline{B}_R$. Then as $\lambda\rightarrow 0$ we
 have,
 $$
 \begin{array}{ll}
 \eta &\leq - \beta (T-\overline{t})^2\\ \eta &< 0,
 \end{array}
 $$
 which yields a contradiction and we have $t_\epsilon \in [0,T)$.\\
 Now we claim that
 \begin{equation}\label{ineg}u'(t_\epsilon,x_\epsilon)-e^{t_\epsilon} h(t_\epsilon,x_\epsilon)>0.
  \end{equation}If not, there exist a subsequence such that $   u'(t_\epsilon,x_\epsilon)-e^{t_\epsilon} h(t_\epsilon,x_\epsilon)\leq 0$, then as $\lambda\rightarrow 0$ we
 have, $u'(\overline{t},\overline{x})-e^{\overline{t}}h(\overline{t},\overline{x})\leq 0$ but from the assumption $u'(\overline{t},\overline{x})-v'(\overline{t},\overline{x})>0$, we deduce that $0\geq u'(\overline{t},\overline{x})-e^{\overline{t}}h(\overline{t},\overline{x})>v'(\overline{t},\overline{x})-h(\overline{t},\overline{x})$. Therefore we have $v'(\overline{t},\overline{x})-e^{\overline{t}}h(\overline{t},\overline{x})<0$, which leads to a contradiction with (\ref{tra}).
 Next let us denote
 \begin{equation*}
 \psi_{\epsilon}(t,x,y)=\frac{1}{2\epsilon}|x-y|^{4}
 +\theta( |x-\overline{x}|^{4}+|y-\overline{x}|^{4})+\beta
 (t-\overline{t})^2.
 \end{equation*}
 Then we have:
\be \left\{
 \begin{array}{lllll}\label{derive}
 D_{t}\psi_{\epsilon}(t,x,y)=2\beta(t-\overline{t}),\\
 D_{x}\psi_{\epsilon}(t,x,y)= \frac{2}{\epsilon}(x-y)|x-y|^{2} +4\theta
 (x-\overline{x})|x-\overline{x}|^{2}, \\
 D_{y}\psi_{\epsilon}(t,x,y)= -\frac{2}{\epsilon}(x-y)|x-y|^{2} +
 4\theta(y-\overline{x})|y-\overline{x}|^{2},\\\\
 B(t,x,y)=D_{x,y}^{2}\psi_{\epsilon}(t,x,y)=\frac{1}{\epsilon}
 \begin{pmatrix}
 a_1(x,y)&-a_1(x,y) \\
 -a_1(x,y)&a_1(x,y)
 \end{pmatrix}+ \begin{pmatrix}
 a_2(x)&0 \\
 0&a_2(y)
 \end{pmatrix} \\\\
 \end{array}
 \right.\ee
 with  $a_1(x,y)=2|x-y|^{2}I+4(x-y)(x-y)^*$  and
 $a_2(x)=4\theta|x-\overline{x}|^{2}I+8xx^*$.\\
 Taking into account (\ref{ineg}) then applying
 the result by Crandall et al. (Theorem 8.3, \cite{MC}) to the
 function
 $$
 (1-\lambda)u'(t,x)-v'(t,y)-\psi_{\epsilon}(t,x,y) $$ at
 the point $(t_\epsilon,x_\epsilon,y_\epsilon)$, for any $\epsilon_1 >0$, we can find
 $c,c_1 \in \R$ and $X,Y \in S(d)$, such that:

 \be \label{lemmeishii}
 \left\{
 \begin{array}{lllll}
 (c,\frac{2}{\epsilon}(x_\epsilon-y_\epsilon)|x_\epsilon-y_\epsilon|^{2} +4\theta
 (x_\epsilon-\overline{x})|x_\epsilon-\overline{x}|^{2},X)
 \in J^{2,+}((1-\lambda)u'(t_\epsilon,x_\epsilon)),\\
 (-c_1,\frac{2}{\epsilon}(x_\epsilon-y_\epsilon)|x_\epsilon-y_\epsilon|^{2} -
 4\theta(y_\epsilon-\overline{x})|y_\epsilon-\overline{x}|^{2},Y)\in J^{2,-}
 (v'(t_\epsilon,y_\epsilon)),\\
 c+c_1=D_{t}\psi_{\epsilon}(t_\epsilon,x_\epsilon,y_\epsilon)=2\beta(t_\epsilon-\overline{t}) \mbox{ and finally }\\
 -(\frac{1}{\epsilon_1}+||B(t_\epsilon,x_\epsilon,y_\epsilon)||)I\leq
 \begin{pmatrix}
 X&0 \\
 0&-Y
 \end{pmatrix}\leq B(t_\epsilon,x_\epsilon,y_\epsilon)+\epsilon_1 B(t_\epsilon,x_\epsilon,y_\epsilon)^2.
 \end{array}\right.\ee\\
Taking now into account (\ref{ineg}), and the
definition of viscosity solution, we get:
$$\begin{array}{l}-c-\frac{1}{2}Tr[\sigma^{*}(t_\epsilon,x_\epsilon)X\sigma(t_\epsilon,x_\epsilon)]-\langle
\frac{2}{\epsilon}(x_\epsilon-y_\epsilon)|x_\epsilon-y_\epsilon|^{2} +4\theta
(x_\epsilon-\overline{x})|x_\epsilon-\overline{x}|^{2},b(t_\epsilon,x_\epsilon)\rangle
+\\\qquad\quad (1-\lambda)u'(t_\epsilon,x_\epsilon)-(1-\lambda)e^{t_\epsilon}f(t_\epsilon,x_\epsilon,e^{-t}u'(t_\epsilon,x_\epsilon),\sigma(t_\epsilon,x_\epsilon)\nabla (e^{-t_\epsilon}u'(t_\epsilon,x_\epsilon)))\leq 0\\
\mbox{ and
}\\c_1-\frac{1}{2}Tr[\sigma^{*}(t_\epsilon,y_\epsilon)Y\sigma(t_\epsilon,y_\epsilon)]-\langle\frac{2}{\epsilon}(x_\epsilon-y_\epsilon)|x_\epsilon-y_\epsilon|^{2} -
4\theta(y_\epsilon-\overline{x})|y_\epsilon-\overline{x}|^{2},b(t_\epsilon,y_\epsilon)\rangle+\\\qquad\quad v'(t_\epsilon,y_\epsilon)-e^{t_\epsilon}f(t_\epsilon,y_\epsilon,e^{-t}v'(t_\epsilon,y_\epsilon),\sigma(t_\epsilon,y_\epsilon)\nabla (e^{-t_\epsilon}v'(t_\epsilon,y_\epsilon)))\geq
0\end{array}$$ which implies that:
\begin{eqnarray}\label{viscder}
\nonumber
&&(1-\lambda)u'(t_\epsilon,x_\epsilon)-v'(t_\epsilon,y_\epsilon)-c-c_1\leq \frac{1}{2}Tr[\sigma^{*}(t_\epsilon,x_\epsilon)X\sigma(t_\epsilon,x_\epsilon)-\sigma^{*}(t_\epsilon,y_\epsilon)Y\sigma(t_\epsilon,y_\epsilon)]\\\nonumber
&&\qquad\qquad\qquad\qquad\qquad\qquad+
\langle\frac{2}{\epsilon}(x_\epsilon-y_\epsilon)|x_\epsilon-y_\epsilon|^{2},b(t_\epsilon,x_\epsilon)-b(t_\epsilon,y_\epsilon)\rangle\\\nonumber
&&\qquad\qquad\qquad\qquad\qquad\qquad+\langle
4\theta
(x_\epsilon-\overline{x})|x_\epsilon-\overline{x}|^{2},b(t_\epsilon,x_\epsilon)\rangle +\langle
4\theta(y_\epsilon-\overline{x})|y_\epsilon-\overline{x}|^{2},b(t_\epsilon,y_\epsilon)\rangle
\\\nonumber
&&\qquad\qquad\qquad\qquad\qquad\qquad+(1-\lambda)e^{t_\epsilon}f(t_\epsilon,x_\epsilon,e^{-t}u'(t_\epsilon,x_\epsilon),\sigma(t_\epsilon,x_\epsilon)\nabla (e^{-t_\epsilon}u'(t_\epsilon,x_\epsilon))\\\nonumber
&&\qquad\qquad\qquad\qquad\qquad\qquad-e^{t_\epsilon}f(t_\epsilon,y_\epsilon,e^{-t}v'(t_\epsilon,y_\epsilon),\sigma(t_\epsilon,y_\epsilon)\nabla (e^{-t_\epsilon}v'(t_\epsilon,y_\epsilon)).\\
&&
\end{eqnarray}

But from (\ref{derive}) there exist two constants $C$ and $C_1$ such
that:
$$||a_1(x_\epsilon,y_\epsilon)||\leq C|x_\epsilon- y_\epsilon|^{2} \mbox{ and }(||a_2(x_\epsilon)||\vee ||a_2(y_\epsilon)||)\leq C_1 \theta.$$
As
$$B= B(t_\epsilon,x_\epsilon,y_\epsilon)= \frac{1}{\epsilon}
\begin{pmatrix}
a_1(x_\epsilon,y_\epsilon)&-a_1(x_\epsilon,y_\epsilon) \\
-a_1(x_\epsilon,y_\epsilon)&a_1(x_\epsilon,y_\epsilon)
\end{pmatrix}+ \begin{pmatrix}
a_2(x_\epsilon)&0 \\
0&a_2(y_\epsilon)
\end{pmatrix}$$
then
$$B\leq \frac{C}{\epsilon}|x_\epsilon - y_\epsilon|^{2}
\begin{pmatrix}
I&-I \\
-I&I
\end{pmatrix}+ C_1 \theta \begin{pmatrix}
I&0 \\
0&I
\end{pmatrix}.$$
It follows that:
\begin{equation*}
B+\epsilon_1 B^2 \leq C(\frac{1}{\epsilon}|x_\epsilon - y_\epsilon|^{2}+
\frac{\epsilon_1}{\epsilon^2}|x_\epsilon - y_\epsilon|^{4})\begin{pmatrix}
I&-I \\
-I&I
\end{pmatrix}+ C_1\theta \begin{pmatrix}
I&0 \\
0&I
\end{pmatrix}
\end{equation*}
where $C$ and $C_1$ which hereafter may change from line to line.
Choosing now $\epsilon_1=\epsilon$, yields the relation
\begin{equation}
\label{ineg_matreciel}
B+\epsilon_1 B^2 \leq \frac{C}{\epsilon}(|x_\epsilon - y_\epsilon|^{2}+|x_\epsilon - y_\epsilon|^{4})\begin{pmatrix}
I&-I \\
-I&I
\end{pmatrix}+ C_1\theta \begin{pmatrix}
I&0 \\
0&I
\end{pmatrix}.
\end{equation}
Now, from the Lipschitz continuity of $\sigma$, (\ref{lemmeishii}) and
(\ref{ineg_matreciel}) we get:
$$\frac{1}{2}Tr[\sigma^{*}(t_\epsilon,x_\epsilon)X\sigma(t_\epsilon,x_\epsilon)-\sigma^{*}(t_\epsilon,y_\epsilon)
Y\sigma(t_\epsilon,y_\epsilon)]\leq \frac{C}{\epsilon}(|x_\epsilon - y_\epsilon|^{4}+|x_\epsilon - y_\epsilon|^{6}) +C_1 \theta.$$ Next by plugging into \eqref{viscder}  we obtain:
\begin{eqnarray*}
&&(1-\lambda)u'(t_\epsilon,x_\epsilon)-v'(t_\epsilon,y_\epsilon)-2\beta(t_\epsilon-\overline{t})\\
&& \qquad\qquad\qquad\qquad\qquad \leq
(1-\lambda)e^{t_\epsilon}f(t_\epsilon,x_\epsilon,e^{-t}u'(t_\epsilon,x_\epsilon),\sigma(t_\epsilon,x_\epsilon)\nabla (e^{-t_\epsilon}u'(t_\epsilon,x_\epsilon))\\
&&\qquad\qquad\qquad\qquad\qquad-e^{t_\epsilon}f(t_\epsilon,y_\epsilon,e^{-t}v'(t_\epsilon,y_\epsilon),\sigma(t_\epsilon,y_\epsilon)\nabla (e^{-t_\epsilon}v'(t_\epsilon,y_\epsilon))\\
&&\qquad\qquad\qquad\qquad\qquad+\frac{C}{\epsilon}(|x_\epsilon - y_\epsilon|^{4}+|x_\epsilon - y_\epsilon|^{6}) +C_1 \theta.
\end{eqnarray*}
By
sending $\epsilon\rightarrow0$, $\lambda\rightarrow0$, $\theta
\rightarrow0$ and taking into account of the continuity of
$f$, we obtain $\eta < 0$ which is
a contradiction.

Now we have $u\leq u\vee h\leq v\wedge h'\leq v.$ which means that if we have $u$ and $\check{u}$ two solutions of \eqref{kkk} then $u\leq \check{u}$ and $\check{u}\leq u$. Hence, obviously, $u=\check{u}$.

\end{document}